\documentclass[11pt]{amsart}

\usepackage{amssymb,amscd,amsthm,amsxtra}
\usepackage{latexsym}
\usepackage{epsfig}

\newtheorem{thm}{Theorem}[section]
\newtheorem{cor}[thm]{Corollary}
\newtheorem{lem}[thm]{Lemma}
\newtheorem{prop}[thm]{Proposition}
\theoremstyle{definition}
\newcommand{\comment}[1]{}

\theoremstyle{remark}
\newtheorem{rem}[thm]{Remark}
\numberwithin{equation}{section}

\begin{document}

\title[Global Well-Posedness for a periodic nonlinear Schr\"odinger equation]
{\bf Global Well-Posedness for a periodic nonlinear Schr\"odinger equation in 1D and 2D}

\author{Daniela De Silva} 
\address{Department of Mathematics, Johns Hopkins University, Baltimore, MD 21218}
\email{\tt  desilva@math.jhu.edu}

\author{Nata\v{s}a Pavlovi\'{c}}
\address{Department of Mathematics, Princeton University, Princeton, NJ 08544-1000}
\email{\tt natasa@math.princeton.edu}

\author{Gigliola Staffilani} 
\address{Department of Mathematics, Massachusetts Institute of Technology, 
Cambridge, MA 02139-4307}
\email{\tt gigliola@math.mit.edu}

\author{Nikolaos Tzirakis}
\address{Department of Mathematics, University of Toronto, Toronto, Ontario, Canada M5S 2E4}
\email{tzirakis@math.toronto.edu}

\date{February 16, 2006} 


\begin{abstract}
The initial value problem for the  $L^{2}$ critical semilinear
Schr\"odinger equation with periodic boundary data is considered.
We show that the problem is globally well posed in $H^{s}({\Bbb
T^{d} })$, for $s>4/9$ and $s>2/3$ in 1D and 2D respectively, 
confirming in 2D a statement of Bourgain in \cite{bo2}. 
We use the ``$I$-method''. This method allows one to introduce a
modification of the energy functional that is well defined for
initial data below the $H^{1}({\Bbb T^{d} })$ threshold. The main
ingredient in the proof is a ``refinement" of the Strichartz's
estimates that hold true for solutions defined on the rescaled space, 
$\Bbb T^{d}_{\lambda} = \Bbb R^{d}/{\lambda \Bbb Z^{d}}$, $d=1,2$.
\end{abstract}
\maketitle

\section{Introduction}

In this paper we study the $L^{2}$ critical Cauchy problem

\begin{align}
&iu_{t}+ \Delta u -|u|^{\frac{4}{d}}u=0, \; \; x \in {\mathbb T}^{d}, \; t \geq 0 \label{ivp1}\\
&u(x,0)=u_{0}(x)\in H^{s}({\Bbb T^{d} }),\label{bc1}
\end{align}
where $d = 1,2$ and $\mathbb T^{d}=\mathbb R^{d}/ \Bbb Z^{d}$ is the d-dimensional torus.

We say that a Cauchy problem is locally well-posed in $H^s$ if for
any choice of initial data $u_0 \in H^s$, there exists a positive
time $T = T(\|u_0\|_{H^s})$ depending only on the norm of the
initial data, such that a solution to the initial value problem
exists on the time interval $[0,T]$, is unique and the solution
map from $H^s_{x}$ to $C^{0}_{t}H^s_{x}$ depends continuously on
the initial data on the time interval $[0,T]$. If $T = \infty$ we
say that a Cauchy problem is globally well-posed.

In the case when $x \in \Bbb R^{d}$ local well-posedness for
\eqref{ivp1}-\eqref{bc1} in $H^{s}({\mathbb R}^{d})$ has been
studied extensively (see, for example, \cite{cw,gv,tk1}). In
particular if one solves the equivalent integral equation by
Picard's fixed point method and controls the nonlinearity in the
iteration process by using Strichartz's type inequalities, then
the problem can be shown to be locally well-posed for all $s>0$.

Bourgain \cite{bo1993} adjusted this approach to the periodic case,
where there are certain difficulties due mainly to a ``lack of
dispersion''. In \cite{bo1993} number theoretic methods were used
to show that \eqref{ivp1}-\eqref{bc1} is locally well-posed in $
H^{s}({\Bbb T^{d} })$, $d=1,2$ for any  $s>0$.

Assuming local existence there are many issues to be addressed
about the behavior of the solution as $t \rightarrow \infty$,
\begin{itemize}
\item global well-posedness/blow-up behavior
\item asymptotic stability
\item behavior of higher order Sobolev norms of smooth solutions.
\end{itemize}

In this note we address the question of global well-posedness for
\eqref{ivp1}-\eqref{bc1}. We recall that solutions of \eqref{ivp1}
satisfy mass conservation
$$ \|u(t)\|_{L^{2}}=\|u_{0}\|_{L^{2}}$$
and smooth solutions also satisfy energy conservation
$$ E(u)(t)=\frac{1}{2}\int |\nabla u(t)|^{2}dx+\frac{d}{2(d+2)}\int |u(t)|^{2+\frac{4}{d}}dx=E(u_{0}).$$
For initial data in $H^{s}(\Bbb T^{d})$ , $s \geq 1$, the energy
conservation together with local well-posedness imply global
well-posedness for $s \geq 1$. However it is a more subtle problem
to extend the global theory to infinite energy initial data.
Bourgain  \cite{bo2} established global well posedness
for \eqref{ivp1} in $H^{s}({\mathbb T})$ for any $s>1/2-$, by
combining a ``normal form" reduction method (see also \cite{bo2004}), the ``$I$-method",
and a refined trilinear Strichartz type inequality. The normal
form reduction is achieved by symplectic transformations that, in
some sense, reduce the nonlinear part of the equation to its
``essential part''. The $I$-method, introduced by J.
Colliander, M. Keel, G. Staffilani, H. Takaoka and T. Tao in
\cite{ckstt1,ckstt4,ckstt5}, is based on the almost conservation
of certain modified Hamiltonians. These two methods together yield
global well-posedness in $H^{s}(\mathbb T)$ with $s\geq 1/2$. The
refined trilinear Strichartz inequality established in \cite{bo2}
is a qualitative estimate needed to achieve global well-posedness
in $H^{s}({\mathbb T})$, $s^{*} < s < 1/2$ for some $s^{*}$.

In this paper we continue to fill in the gap between what is known
locally and what is known globally in $H^{s}(\Bbb T^{d})$ for
$s>0$, when $d=1,2$. Our approach is based on an implementation of
the $I$-method itself adjusted to the periodic setting via
elementary number theoretic techniques.

In order to present our method, we briefly review global well-posedness results on ${\mathbb R}^d$.
 Using an approximation of the modified energy in the $I$-method, Tzirakis \cite{tz} showed that the Cauchy problem
 \eqref{ivp1}-\eqref{bc1} is globally well posed in $H^{s}({\Bbb R})$ for any $s>4/9$.
 For $u_{0} \in H^{s}({\Bbb R^{2}})$ the best known global well posedness result for
\eqref{ivp1}-\eqref{bc1} is in \cite{ckstt5} where the authors proved global well posedness
for any $s>4/7$. In both results, apart from the application of the $I$-method,
a key element in the proof is the existence of a bilinear refined Strichartz's
estimate due to Bourgain, \cite{bo3} (see also \cite{cdks}
\cite{ckstt1}). In general dimensions $d \geq 2$ this estimate reads as follows. Let
 $f$ and $g$ be any two Schwartz functions whose Fourier transforms are supported in $|k| \sim N_{1}$ and $|k| \sim N_{2}$
respectively. Then we have
$$\|U_{t}f U_{t}g\|_{L_{t}^{2}L_{x}^{2}( \Bbb R\times \Bbb R^{d})} \leq C
\frac{N_{2}^{\frac{d-1}{2}}}{N_{1}^{\frac{1}{2}}}\|f\|_{L^{2}( \Bbb
R^{d})}\|g\|_{L^{2}( \Bbb R^{d})},$$
where $U_t$ denotes the solution operator associated to the linear
Schr\"{o}dinger equation. In one dimension the estimate fails if
the two frequencies are comparable but continues to hold if the
frequencies are separated ($N_1 >> N_2$). Such an estimate is very
useful when $f$ is in high frequency and $g$ is in low frequency
since we can move derivatives freely to the low frequency factor.

The main difficulty in obtaining global well-posedness results in the periodic context
is exactly the absence of a quantitative refinement of the bilinear Strichartz's estimate.
The reader can consult the paper of Kenig, Ponce and Vega, \cite{kpv}, where the difference between the real
and the periodic case is clearly exposed when one tries to prove bilinear estimates in different functional spaces.

In order to overcome the
non-availability of a quantitative refined Strichartz's estimate,
we develop a different strategy, closer in spirit to the approach
in \cite{ckstt4}. By exploiting the scaling symmetry of the
equation, we analyze \eqref{ivp1} on $\Bbb T_{\lambda}^{d} \times \Bbb
R$, $d=1,2,$ where $\Bbb T_{\lambda}^{d}={\Bbb R}^{d} /\lambda {\Bbb Z}^{d}$. The main
novelty in this analysis consists in establishing a bilinear
Strichartz's inequality for $\lambda$-periodic functions which are
well separated in frequency space. The constant in the inequality
is quantified in terms of $\lambda$. As $\lambda \rightarrow
\infty$, our estimate reduces to the refined bilinear Strichartz's
inequality\footnote{ In a way we are introducing ``more
dispersion'' by rescaling the problem.} on ${\mathbb R}^{d}$, $d=1,2$, see, for example, 
\cite{bo3, ckstt1, ot1998}. Such an estimate allows us to use 
the $I$-method machinery in an efficient way.

More precisely, when $d=1$,  beside rescaling, we follow the argument 
in \cite{tz} (see also \cite{ckstt3}), where one applies  
the $I$-operator to the equation on ${\Bbb T}_{\lambda}$, 
and defines a modified second energy functional as the
energy corresponding to the  new ``$I$-system''. We prove that
such modified second energy is ``almost conserved'', that is the time
derivative of this new energy decays with respect to a very large
parameter $N$. Roughly speaking, $N$ denotes the stage at which
the $I$-operator stops behaving like the identity and starts
smoothing out the solution. We prove a local existence result for
the ``$I$-system'', under a smallness assumption for the initial
data, which is guaranteed by choosing $\lambda$ in terms of $N$.
The decay of the modified energy enables us to iterate the local
existence preserving the same bound for $\|u^{\lambda}\|_{H^{s}}$
during the iteration process. By undoing the scaling we obtain 
polynomial in time bounds for $u$ in
${\Bbb T} \times {\Bbb R}^{+}$,  and this immediately implies global well
posedness.

The precise statement of  our 1D result reads as follows.
\begin{thm}\label{main1D}
The initial value problem \eqref{ivp1}-\eqref{bc1} is globally well-posed in $H^{s}(\mathbb T)$ for
$s > \frac{4}{9}$.
\end{thm}

For the two dimensional case  Bourgain already announced in \cite{bo2}
that the $I$-method based only on the first energy would give global
well-posedness for $s>2/3$. While we were explicitly writing up the
calculations to recover this claim,  we noticed  that in one particular 
case\footnote{See Case IIIb) in the proof of Proposition \ref{energybound}.}
of the estimate of the first energy, a better Strichartz inequality
was needed to successfully conclude the argument. We  then proceeded by determining a
qualitative ``$\epsilon$-refined" Strichartz type estimate, see Proposition \ref{2Dcounting}, 
which allowed us to implement the $I$-method described above to obtain
indeed  the following result:

\begin{thm}\label{main2D}
The initial value problem \eqref{ivp1}-\eqref{bc1} is globally
well-posed in $H^{s}({\mathbb T}^{2})$ for $s > \frac{2}{3}$.
\end{thm}

We remark that,  as we mentioned above,  in 1D we introduce an approximation to the modified
energy, by adding correction terms, in the spirit of \cite{tz} and
\cite{ckstt3}. In
the 2D problem we did not use correction terms, since the Fourier
multipliers corresponding to the approximated modified energy
would be singular.  This singularity is not only caused by the
presence of zero frequencies, but also  by  orthogonality issues. 
This is a whole new ground that we are exploring in a different paper.

\subsection*{Organization of the paper} In section 2 we review the notation and
some known Strichartz estimates for the periodic Schrodinger equations
in 1D and 2D. In section 3 we present the proof of Theorem \ref{main1D},
while in section 4 we give a proof of Theorem \ref{main2D}.

\subsection*{Acknowledgments} We would like to thank Akshay Venkatesh
and Terence Tao for valuable discussions. N.P. received partial support from NSF grant 
DMS-0304594. N.T. received partial support from NSF grant DMS-0111298. Portions of this work
were carried out at the Institute for Advanced Study (IAS) in Princeton and at the  
Mathematical Sciences Research Institute (MSRI) at Berkeley. The authors would like to thank
these institutions for their hospitality.

\section{Notation}

In what follows we use $A \lesssim B$ to denote an estimate of the form $A\leq CB$ for some constant $C$.
If $A \lesssim B$ and $B \lesssim A$ we say that $A \sim B$. We write $A \ll B$ to denote
an estimate of the form $A \leq cB$ for some small constant $c>0$. In addition $\langle a \rangle:=1+|a|$ and
$a\pm:=a\pm \epsilon$ with $0 < \epsilon <<1$.

%We will work with $\lambda-$periodic functions, thus we define
%$(dk)_{\lambda}$ to be the normalized counting measure on
%$\frac{1}{\lambda}\Bbb Z^2$:
%$$\int a(x)(dk)_{\lambda}=\frac{1}{\lambda^2}\sum_{k \in \frac{1}{\lambda}\Bbb Z}a(k).$$
%Then define the Fourier transform of $f(x) \in L_{x\in
%[0,\lambda]}^{1}$ by
%$$\hat{f}(k)=\int_{0}^{\lambda}e^{-2\pi ikx}f(x)dx$$
%For an appropriate class of functions we know by the Fourier
%inversion formula that:
%$$f(x)=\int e^{2\pi ikx}\hat{f}(k)(dk)_{\lambda}.$$

We recall that the equation \eqref{ivp1} is $L^{2}$ invariant under the following scaling
$(x,t)\rightarrow \frac{1}{\lambda^{\frac{d}{2}}}(\frac{x}{\lambda},\frac{t}{\lambda^{2}})$.
 Thus if $u(x,t)$ solves \eqref{ivp1} on $\Bbb T^{d} \times \mathbb{R}$
then
$$u^{\lambda}(x,t)=\frac{1}{\lambda^{\frac{d}{2}}}u(\frac{x}{\lambda},\frac{t}{\lambda^{2}})$$
is a solution of (\ref{ivp1}) in ${\Bbb T^{d}}_{\lambda} \times \Bbb R$
where $\Bbb T^{d}_{\lambda} = \Bbb R^{d}/{\lambda \Bbb Z^{d}}.$

Since in our argument we exploit the scaling symmetry let us recall some properties of
$\lambda-$periodic functions. Define $(dk)_{\lambda}$ to be the normalized counting
measure on $(\frac{1}{\lambda}{\Bbb Z})^{d}$:
$$\int a(k)(dk)_{\lambda}=\frac{1}{{\lambda}^{d}}\sum_{k \in (\frac{1}{\lambda}{\Bbb Z})^{d}}a(k).$$
We define the Fourier transform of $f(x) \in L_{x\in [0,\lambda]^{d}}^{1}$ by
$$\hat{f}(k)=\int_{[0,\lambda]^{d}} e^{-2\pi ikx}f(x)dx.$$
For an appropriate class of functions the following Fourier inversion formula holds:
$$f(x)=\int e^{2\pi ikx}\hat{f}(k)(dk)_{\lambda}.$$
Moreover we know that the following identities are true:\\
\begin{enumerate}
\item $\|f\|_{L^{2}([0,\lambda]^{d})}=\|\hat{f}\|_{L^{2}((dk)_{\lambda})}$, (Plancherel)\\
\item $\int_{[0,\lambda]^d} f(x)\bar{g}(x)dx=\int \hat{f}(k)\bar{\hat{g}}(k) (dk)_{\lambda}$, (Parseval)\\
\item $\widehat{fg}(k)=\hat{f}\star _{\lambda} \hat{g}(k)=\int \hat f(k-k_{1})\hat g(k_{1})(dk_{1})_{\lambda}$,\\
\end{enumerate}

We define the Sobolev space $H^{s} = H^{s}([0,\lambda]^{d})$ as the space equipped with the norm
$$\|f\|_{H^{s}}=\|\langle k \rangle^{s}\hat{f}(k)\|_{L^{2}((dk)_{\lambda})}.$$
We write $U_{\lambda}(t)$ for the solution operator to the linear Schr\"odinger equation
$$iu_t - \Delta u = 0,\; \; x \in [0,\lambda]^{d},$$
that is
$$U_{\lambda}(t)u_{0}(x)=\int e^{2\pi ikx-(2\pi k)^{2}it}\widehat{u_0}(k)(dk)_{\lambda}.$$
We denote by $X^{s,b} = X^{s,b}(\mathbb{T}^{d}_{\lambda} \times \mathbb R)$ the completion of
${\mathcal S}({\mathbb T}^{d}_{\lambda} \times \mathbb R)$ with respect to the following norm,
see, for example, \cite{g}

$$\|u\|_{X^{s,b}}= \|U_{\lambda}(-t)u\|_{H_{x}^{s}H_{t}^{b}}
=\|\langle k \rangle^{s} \langle \tau -4\pi ^{2}k^{2} \rangle^{b}\tilde{u}(k,\tau)\|_{L_{\tau}^{2}L_{(dk)_{\lambda}}^{2}},$$

where  $\tilde{u}(k,\tau)$ is the space-time Fourier Transform

$$\tilde{u}(k,\tau)=\int \int_{[0,\lambda]^{d}} e^{-2 \pi i(k \cdot x + \tau t)}u(x,t)dxdt.$$

Furthermore for a given time interval $J$, we define
$$\|f\|_{X^{s,b}_{J}}=\inf_{g  = f \mbox{ on } J}  \|g\|_{X^{s,b}}.$$
Often we will drop the subscript $J$.

In what follows we will need the following known estimates. Since
$$\widehat{U_{\lambda}(t)u_{0}}(k)=e^{(2\pi k)^{2}it}\widehat{u_{0}}(k)$$
we have that
$$\|U_{\lambda}(t)u_{0}\|_{L_{t}^{\infty}L_{x}^{2}} \lesssim \|u_{0}\|_{L^{2}}.$$
Hence
\begin{equation} \label{inf2}
\|u\|_{L_{t}^{\infty}L_{x}^{2}}
=\|U_{\lambda}(t)U_{\lambda}(-t)u\|_{L_{t}^{\infty}L_{x}^{2}}
\lesssim \|U_{\lambda}(-t)u\|_{L^{\infty}_t L_{x}^{2}} \lesssim
\|u\|_{X^{0,1/2+}},
\end{equation}
where in the last inequality we applied the definition of the
$X^{s,b}$ spaces, the basic estimate $\|u\|_{L^{\infty}} \leq
\|\hat{u}\|_{L^{1}}$, and the Cauchy-Schwartz inequality.

\subsection*{1D estimates}
If $u$ is on ${\mathbb{T}}^{1}$ then the estimate \eqref{inf2} 
combined with the Sobolev embedding
theorem implies that
\begin{equation} \label{inf}
\|u\|_{L_{t}^{\infty}L_{x}^{\infty}} \lesssim
\|u\|_{X^{1/2+,1/2+}}.
\end{equation}
Also the following linear Strichartz's estimates was obtained in
\cite{bo1993} in the case of the torus:
\begin{equation} \label{L4}
\|u\|_{L_{t}^{4}L_{x}^{4}} \lesssim \|u\|_{X^{0,b}},
\end{equation}
for any $b>\frac{3}{8}$, and 
\begin{equation} \label{1L6}
\|u\|_{L_{t}^{6}L_{x}^{6}} \lesssim \|u\|_{X^{0+,1/2+}}.
\end{equation}
 
We note that \eqref{inf} remains true for the $\lambda$-periodic problem. This is the case
also for \eqref{L4}. The proof is essentially in \cite{gr}. In fact, 
it is enough to show, by the standard $X^{s,b}$ method, that
$$\frac{1}{\lambda}\sup_{(\xi, \tau)\in \frac{1}{\lambda}\Bbb Z \times \Bbb R}\sum_{k_{1}\in \frac{1}{\lambda}\Bbb Z}
\langle \tau+k_{1}^{2}+(k-k_{1})^{2} \rangle^{1-4b}\lesssim C.$$
This is done in \cite{gr}, the only difference being that there are $O(\lambda)$ numbers $k\in \frac{1}{\lambda}\Bbb Z$,
such that $|k+x_{0}|<1$ and $|k-x_{0}|<1$, where $x_{0}^{2}=|2\tau +\xi^{2}|$. But then summing the above
 series, using Cauchy-Schwartz inequality and the fact that $b>3/8$, one gets that the left hand side of the inequality
 is
$$\lesssim \frac{1}{\lambda}(c_{1}\lambda+c_{2}) \lesssim C$$ for any $\lambda>1$. 
So from now on we will use \eqref{L4} for the $\lambda$-periodic solutions without any further comment.
On the other hand, using scaling we can see that for the $\lambda$-periodic solutions, \eqref{1L6} takes the form
\begin{equation} \label{L6}
\|u\|_{L_{t}^{6}L_{x}^{6}} \lesssim
\lambda^{0+}\|u\|_{X^{0+,1/2+}}.
\end{equation}
If we interpolate equations
\eqref{inf} and \eqref{L6} we get
\begin{equation} \label{int1}
\|u\|_{L_{t}^{p}L_{x}^{p}} \lesssim \lambda^{0+}\|u\|_{X^{\alpha_1(p),1/2+}},
\end{equation}
with $\alpha_1(p) = (\frac{1}{2}-\frac{3}{p})+$ and $6 \leq p \leq \infty$.

 \subsection*{2D estimates}
On ${\mathbb{T}}^{2}$ Bourgain proved, \cite{bo1993},
\begin{equation} \label{1L44}
\|u\|_{L_{t}^{4}L_{x}^{4}} \lesssim
\|u\|_{X^{0+,1/2+}}.
\end{equation}
Again using scaling we can prove that for the $\lambda$-periodic solutions, \eqref{1L44} takes the form
\begin{equation} \label{L44}
\|u\|_{L_{t}^{4}L_{x}^{4}} \lesssim \lambda^{0+}
\|u\|_{X^{0+,1/2+}}.
\end{equation}
The estimate \eqref{inf2} together with Sobolev embedding gives:
\begin{equation} \label{energy}
\|u\|_{L_{t}^{\infty}L_{x}^{\infty}} \lesssim \|u\|_{X^{1+,1/2+}}.
\end{equation}
Hence, by interpolation, we get
\begin{equation} \label{interpol}
\|u\|_{L_{t}^{p}L_{x}^{p}} \lesssim \lambda^{0+}
\|u\|_{X^{\alpha_2(p),1/2+}},
\end{equation}
with $\alpha_{2}(p)=(1-\frac{4}{p})+$ and $4 \leq p \leq \infty.$

\begin{rem} \label{decrm} 
({\bf{Decomposition remark}})
Our approach to prove Theorem \ref{main1D} and Theorem \ref{main2D} is based on obtaining certain multilinear estimates in appropriate functional spaces which are $L^2$-based. Hence, whenever we perform a Littlewood-Paley decomposition of a function
we shall assume that the Fourier transforms of the Littlewood-Paley pieces are positive. Moreover,
we will ignore the presence of conjugates.
\end{rem}

\section{The I-method and the proof of Theorem \ref{main1D}}

In this section we present the proof of Theorem \ref{main1D}. We
start by recalling the definition of the operator $I$ introduced
by Colliander et al, see, for
example, \cite{ckstt1,ckstt5}. For $s<1$ and a parameter $N >>1$
let $m(k)$ be the restriction to $\Bbb Z$ of the following smooth
monotone multiplier:
\[m(k):= \left\{\begin{array}{ll}
1 & \mbox{if $|k|<N$}\\
(\frac{|k|}{N})^{s-1} & \mbox{if $|k|>2N$}
\end{array}
\right.\] We define the multiplier operator $I:H^{s} \rightarrow
H^{1}$ by
$$\widehat{Iu}(k)=m(k)\hat{u}(k).$$
The operator $I$ is smoothing of order $1-s$ and we have that:
\begin{equation}
\|u\|_{X^{s_{0},b_{0}}} \lesssim \|Iu\|_{X^{s_{0}+1-s,b_{0}}}
\lesssim N^{1-s}\|u\|_{X^{s_{0},b_{0}}}
\end{equation}
for any $s_{0},b_{0}\in {\mathbb R}$.

%$\bf{ Remark.}$ It is shown in \cite{ckstt2} that if
%$$\|uv\|_{X^{s,b-1}} \lesssim \|u\|_{X^{s,b}} \|v\|_{X^{s,b}}$$
%then
%$$\|I(uv)\|_{X^{1,b-1}} \lesssim  \|Iu\|_{X^{1,b}} \|Iv\|_{X^{1,b}}$$
%where the constants in the inequality above are independent of $N$. From now on we use this fact
%and refer to it as the ``interpolation lemma''. For details see \cite{ckstt2}.

For $u \in H^s$ we set
\begin{equation}\label{first}E^{1}(u)=E(Iu),\end{equation} where 
$$E(u)(t)=\frac{1}{2}\int |\nabla u(t)|^{2}dx+\frac{1}{6}\int |u(t)|^{6}dx=E(u_{0}).$$ We refer
to $E^{1}(u)$ as the first modified energy. As observed by
Colliander et al a hierarchy of
modified energies can be formally considered for different nonlinear
dispersive equations. The goal of the $I$-method is to prove that
the modified energies are ``almost conserved" i.e. they decay in
time with respect to $N$. Since in 1D we base our approach on the
analysis of a second modified energy, it is appropriate to collect
some facts concerning the calculus of multilinear forms used to
define the hierarchy, see, for example \cite{tz}.

If $n\geq2$ is an even integer we define a spatial multiplier of
order $n$ to be the function $M_{n}(k_{1},k_{2},\ldots ,k_{n})$ on
$\Gamma _{n}=\{(k_{1},k_{2},\ldots ,k_{n} ) \in
\frac{1}{\lambda}\Bbb Z^{n}:k_{1}+k_{2}+ \ldots +k_{n}=0\}$ which
we endow with the standard measure $\delta(k_{1}+k_{2}+ \ldots
+k_{n})$. If $M_{n}$ is a multiplier of order $n$, $1 \leq j \leq
n$ is an index and $l \geq 1$ is an even integer we define the
elongation $X_{j}^{l}(M_{n})$ of $M_{n}$ to be the multiplier of
order $n+l$ given by
$$X_{j}^{l}(M_{n})(k_{1},k_{2},\ldots ,k_{n+l})=
M_{n}(k_{1},\ldots,k_{j-1},k_{j}+\ldots+k_{j+l},k_{j+l+1}, \ldots,
k_{n+l}).$$ In addition if $M_{n}$ is a multiplier of order $n$
and $f_{1},f_{2},...,f_{n}$ are functions on ${{\mathbb
T}_{\lambda}}$ we define
$$\Lambda_{n}(M_{n};f_{1},f_{2},...,f_{n})=
\int_{\Gamma_{n}}M_{n}(k_{1},k_{2},\ldots ,k_{n})\prod_{i=1}^{n}\hat f_{j}(k_{j}),$$
where we adopt the notation
$\Lambda_{n}(M_{n};f)=\Lambda_{n}(M_{n};f,\bar f,...,f,\bar f)$.
Observe that  $\Lambda_{n}(M_{n};f)$ is invariant under
permutations of the even $k_{j}$ indices, or of the odd $k_{j}$
indices.

If $f$ is a solution of \eqref{ivp1} the following differentiation
law holds for the multilinear forms $\Lambda_{n}(M_{n};f)$:
\begin{equation}\label{diff}
\partial_{t}\Lambda_{n}(M_{n})=i\Lambda_{n}(M_{n}\sum_{j=1}^{n}(-1)^{j}k_{j}^{2})+i\Lambda_{n+4}
(\sum_{j=1}^{n}(-1)^{j}X_{j}^{4}(M_{n})).
\end{equation}
Note that in this notation the first modified energy \eqref{first}
reads as:
$$E^1(u)=\frac{1}{2}\int |\partial_{x}Iu|^{2}dx+\frac{1}{6}\int |Iu|^{6}dx=
-\frac{1}{2}\Lambda_{2}(m_{1}k_{1}m_{2}k_{2})+\frac{1}{6}
\Lambda_{6}(m_{1}...m_{6})$$ where $m_{j}=m(k_{j})$.

We define the second modified energy
$$E^{2}(u)=-\frac{1}{2}\Lambda_{2}(m_{1}k_{1}m_{2}k_{2})
+\frac{1}{6}\Lambda_{6}(M_{6}(k_{1},k_{2},...,k_{6})),$$
where $M_{6}(k_{1},k_{2},...,k_{6})$ is the following multiplier:
\begin{equation}\label{M6}
M_{6}(k_{1},k_{2},...,k_{6})=
\frac{m_{1}^{2}k_{1}^{2}-m_{2}^{2}k_{2}^{2}+
m_{3}^{2}k_{3}^{2}-m_{4}^{2}k_{4}^{2}+m_{5}^{2}k_{5}^{2}
-m_{6}^{2}k_{6}^{2}}{k_{1}^{2}-k_{2}^{2}+k_{3}^{2}-k_{4}^{2}+k_{5}^{2}-k_{6}^{2}}.
\end{equation}
Notice that the zero set of the denominator corresponds to the
resonant set for a six-waves interaction. We remark that $M_6$
contains more ``cancellations" than the multiplier $m_1...m_6$
that appears in $E^1$.

The differentiation rule \eqref{diff} together with the
fundamental theorem of calculus implies the following Lemma, which
will be used to prove that $E^2$ is almost conserved.

\begin{lem}\label{fundcal}
Let $u$ be an $H^{1}$ solution to $(\ref{ivp1})$. Then for any
$T\in \bf{R}$ and $\delta >0$ we have
\begin{equation} \label{increment}
E^{2}(u(T+\delta))-E^{2}(u(T))=\int_{T}^{T+\delta}\Lambda_{10}(M_{10};u(t))dt,
\end{equation}
with $M_{10}= c \sum \{
M_{6}(k_{abcde},k_{f},k_{g},k_{h},k_{i},k_{j})
-M_{6}(k_{a},k_{bcdef},k_{g},k_{h},k_{i},k_{j})+ \\ \\
M_{6}(k_{a},k_{b},k_{cdefg},k_{h},k_{i},k_{j})-M_{6}(k_{a},k_{b},k_{c},k_{defgh},k_{i},k_{j})+\\ \\
M_{6}(k_{a},k_{b},k_{c},k_{d},k_{efghj},k_{j})-M_{6}(k_{a},k_{b},k_{c},k_{d},k_{e},k_{fghij})\}$,\\ \\
where the summation runs over all permutations
$\{a,c,e,g,i\}=\{1,3,5,7,9\}$ and $\{b,d,f,h,j\}
\\=\{2,4,6,8,10\}$. Furthermore if $|k_{j}|\ll N$ for all $j$ then the multiplier $M_{10}$ vanishes.
\end{lem}

As it was observed in \cite{tz} where the equation \eqref{ivp1}
was considered on $\mathbb R$ one has

\begin{prop}\label{Multb}
The multiplier $M_6$ defined in $(\ref{M6})$ is bounded on its
domain of definition.
\end{prop}
\noindent For the proof see \cite{tz}.

\

Now we proceed to present the steps leading to the proof of
Theorem \ref{main1D}.
\subsection*{Step 1: Local well-posedness for the ${\mathbf{I}}$-system}
The first step towards the proof of Theorem \ref{main1D} is to
apply the $I$-operator to \eqref{ivp1} and prove a local
well-posedness result for the $I$-initial value problem
\begin{align}\begin{split}\label{Iivp}
&iIu_{t}+Iu_{xx}-I(|u|^{4}u)=0\\
&Iu(x,0)=Iu_{0}(x)\in H^{1}({\Bbb T_{\lambda}}),t\in {\Bbb R}.
\end{split}
\end{align}
In order to obtain the local well-posedness we need the following
technical lemma:

\begin{lem}\label{norms}
Let $\eta \in C_{0}^{\infty}$ be a bump function that is supported
on $[-2,2]$ and equals 1 on $[-1,1]$ and set
$\eta_{\delta}(t)=\eta(\frac{t}{\delta})$. For $b,b^{'} \in \Bbb
R$ with $-1/2<b^{'}\leq 0 \leq b \leq b^{'}+1$
 and $\delta \leq 1$ we have:
\begin{enumerate}
\item $\|\eta_{1}(t)U_{\lambda}(t)u_{0}\|_{X^{s,b}} \lesssim \|u_{0}\|_{H^{s}},$
\item $\|\eta_{\delta}(t)\int_{0}^{\delta}U_{\lambda}(t-\tau)F(\tau)d\tau\|_{X^{s,b}}
\lesssim \delta^{1+b^{'}-b} \|F\|_{X^{s,b'}}$.
\end{enumerate}
\end{lem}

\begin{proof}

(1) We compute (see, for example, \cite{g})
\begin{align*}
\|\eta_{1}(t)U_{\lambda}(t)u_{0}\|_{X^{s,b}} & = \|U_{\lambda}(-t)\eta_{1}(t)U_{\lambda}(t)u_{0}\|_{H_{t}^{b}H_{x}^{s}} \\
& = \|\eta_{1}(t)u_{0}\|_{H_{t}^{b}H_{x}^{s}} \\
& = \|u_{0}\|_{H^{s}}\|\eta_{1}(t)\|_{H_{t}^{b}} \\
& \sim \|u_{0}\|_{H^{s}}.
\end{align*}
\\
(2) The following inequality is true and its proof can be found in \cite{g}, \cite{kpv},

$$\|\eta_{\delta}(t)\int_{0}^{\delta}G(\tau)d\tau\|_{H_{t}^{b}} \lesssim \delta^{1+b^{'}-b} \|G\|_{H_{t}^{b^{'}}}.$$
But then $(2)$ follows if we apply the above inequality for fixed $\xi$, multiplying by $<\xi>^{2s}$ and taking
the $L^{2}$ norm in $\xi$.
\end{proof}

\subsection*{Remark} It is shown in \cite{ckstt2} that if
$$\|uv\|_{X^{s,b-1}} \lesssim \|u\|_{X^{s,b}} \|v\|_{X^{s,b}},$$
then
$$\|I(uv)\|_{X^{1,b-1}} \lesssim  \|Iu\|_{X^{1,b}} \|Iv\|_{X^{1,b}},$$
where the constants in the inequality above are independent of
$N$. From now on we use this fact and refer to it as the
``invariant lemma''. For details see Lemma 12.1 in \cite{ckstt2}.
\\
\\
Now we present the local well-posedness result for \eqref{Iivp}. 
Because of the loss in the periodic Strichartz estimates \eqref{L6} 
represented by the $\lambda^{0+}$ factor, the local well-posedness is 
not as straightforward as in the real case. (There one can
essentially use the Leibnitz rule for the fractional derivatives 
and the Cauchy-Schwartz inequality along with interpolated Strichartz's 
estimates to bound the nonlinearity). In 1D we can still follow the same approach
thanks to \eqref{L4}, which is valid for the $\lambda$-periodic problem.
In 2D for the  $\lambda$-periodic
problem we need a more complicated argument that is due to
Bourgain, \cite{bo1}, and uses a local variant of the periodic Strichartz 
inequality that is proved there. For details see \cite{bo1}. As a final comment, 
note that our local well-posedness results are not optimal but
they suffice for the purposes of the global theory that we establish. 
Thus without aiming at sharpness, we
prove local well-posedness for $s>5/16$ in 1D (see Proposition \ref{lwp} and 
Corollary \ref{lwp1}) and $s>7/20$ in 2D (see Proposition \ref{2Dlwp} and
Corollary \ref{n2Dlwp}).

\begin{prop}\label{lwp}
Consider the $I$-initial value problem \eqref{Iivp}, $\lambda=1$. If
$\|Iu_{0}\|_{H^{1}} \lesssim 1$, then $\eqref{Iivp}$ is locally
well-posed for any $s>\frac{5}{16}$ in $[0,\delta]\sim [0,1]$.
\end{prop}

\begin{proof}
By Duhamel's formula and Lemma \ref{norms} we have that
\begin{equation} \label{duh1}
\|Iu\|_{X^{1,1/2+}} \lesssim
\|Iu_{0}\|_{H^{1}}+\delta^{\frac{1}{2}-\epsilon}\|I(|u|^{4}u)\|_{X^{1,-1/2+2\epsilon}}.
\end{equation}
We shall prove that
$$\|I(|u|^{4}u)\|_{X^{1,-1/2+2\epsilon}}\lesssim \|Iu\|_{X^{1,1/2+}}^{5}.$$
By the invariant lemma we know that to prove such an estimate
it suffices to prove
$$\||u|^{4}u\|_{X^{s,-1/2+2\epsilon}}\lesssim \|u\|_{X^{s,1/2+}}^{5}.$$
By H\"older's inequality combined with the Leibnitz rule (\cite{kpvkdv})
we have
$$\||u|^{4}u\|_{X^{s,-1/2+2\epsilon}} \leq \||u|^{4}u\|_{X^{s,0}} \lesssim
\|J^{s}u\|_{L_{t}^{4}L_{x}^{4}}\|u\|^4_{L_{t}^{16}L_{x}^{16}},$$
where $J^{s}$ is the Bessel potential of order $s$. Using
\eqref{int1} with $\lambda = 1$ we obtain
$$\|u\|_{L_{t}^{16}L_{x}^{16}}\lesssim \|u\|_{X^{5/16+,1/2+}} \lesssim \|u\|_{X^{s,1/2+}},$$
for $s>\frac{5}{16}$. Note that we also use \eqref{L4} in order to bound $\|J^{s}u\|_{L_{t}^{4}L_{x}^{4}}$. Thus
$$\||u|^{4}u\|_{X^{s,-1/2+2\epsilon}}\lesssim \|u\|_{X^{s,1/2+}}^{5}$$
which combined with \eqref{duh1} gives
$$\|Iu\|_{X^{1,1/2+}} \lesssim \|Iu_{0}\|_{H^{1}}+\delta^{\frac{1}{2}-\epsilon}\|Iu\|_{X^{1,1/2+}}^{5}.$$
Now by standard nonlinear techniques, see, for example,
\cite{ckstt3} we have that for $\delta \sim 1$
$$\|Iu\|_{X^{1,1/2+}} \lesssim \|Iu_{0}\|_{H^{1}}.$$
\end{proof}
\begin{cor}\label{lwp1}
Consider the $I$-initial value problem \eqref{Iivp}. If
$\|Iu_{0}\|_{H^{1}} \lesssim 1$, then $\eqref{Iivp}$ is locally
well-posed for any $s>\frac{5}{16}$ in $[0,\delta]\sim [0,\frac{1}{\lambda^{a\epsilon}}]$, where $a$ is a fixed positive
 integer ($a \sim 20$ will do it).
\end{cor}
\begin{proof}
The proof is line by line the same as in the previous proposition.
The important step is that in the iteration process where
$x \in [0,\lambda]$, $\delta$ has a fixed power and this fact 
gives us a local well-posedness result on the time interval
$[0, \frac{1}{\lambda^{a\epsilon}}]$, where 
$\frac{1}{\lambda^{a\epsilon}}\sim \frac{1}{\lambda^{0+}}$ 
for all practical purposes.
\end{proof}

\subsection*{Step 2: ${\mathbf{E^{2}(u)}}$ is a small perturbation of ${\mathbf{E^{1}(u)}}$ } Now we shall prove that the
second energy $E^2(u)$ is a small perturbation of the first
energy $E^1(u)$. This is the content of the following proposition.

\begin{prop}\label{edif}
Assume that $u$ solves $(\ref{Iivp})$ with $s \geq 2/5$. Then
$$E^{2}(u) = E^{1}(u)+O(1/N)\|Iu\|_{H^{1}}^{6}.$$
Moreover if $\|Iu\|_{H^{1}}=O(1)$ then
\begin{equation} \label{corder} 
\|\partial_{x}Iu\|_{L^{2}}^{2}\lesssim E^{2}(u) + O(\frac{1}{N}).
\end{equation} 
\end{prop}

\begin{proof}
By the definition of first and second modified energy we have
\begin{equation} \label{E12}
E^{2}(u)=-\frac{1}{2}\Lambda_{2}(m_{1}k_{1}m_{2}k_{2}) + \frac{1}{6}\Lambda_{6}(M_{6})
=E^{1}(u)+ \frac{1}{6}\Lambda_{6}(M_{6}-\prod_{i=1}^{6}m_{i} ).
\end{equation}
Therefore it suffices to prove the following pointwise in time
estimate
\begin{equation}\label{perturbation}
|\Lambda_{6}(M_{6}-\prod_{i=1}^{6}m_{i} )|(t) \lesssim O(1/N)
\|Iu(\cdot, t)\|_{H^{1}_x}^{6}.\end{equation} Combining the
\textbf{Decomposition Remark} \ref{decrm} 
with the fact that $M_6$ is bounded
(by Proposition \ref{Multb}) and that $m$ is bounded (by its
definition) it is enough to show that

\begin{equation}\nonumber 
\int_{\Gamma_{6}} \prod_{j=1}^{6}\widehat{u}(k_{j}) 
\lesssim \frac{1}{N}\|Iu(\cdot,t)\|_{H^{1}_x}^{6}.
\end{equation} 

Towards this aim, we break all
the functions into a sum of dyadic constituents $u_{j}$, each
with frequency support $\langle k \rangle \sim 2^{j}, j=0,...$. We
will be able to sum over all frequency pieces $u_j$ of $u$
since our estimate will decay geometrically in these frequencies.
Hence we need to show the following:
\begin{equation}\label{pieceb}
\int_{\Gamma_{6}} \prod_{j=1}^{6}\widehat{u_j}(k_{j}) \lesssim
\frac{1}{N}(N_1...N_6)^{0-}\prod_{j=1}^{6}\|Iu_j(\cdot,
t)\|_{H^{1}_x}^{6}\end{equation} where $u_j$ is supported around
$\langle k \rangle \sim N_j=2^{h_j}$ for some $h_j$'s.

Denote by $N_{1}^{\star} \geq N_{2}^{\star} \geq N_{3}^{\star}
\geq N_{4}^{\star} \geq N_{5}^{\star} \geq N_{6}^{\star}$ the
decreasing rearrangement of $N_1,...,N_6$, and by $u_j^\star$ the
function $u_j$ supported in Fourier space around $N_j^\star$,
$j=1,...,6.$

Since we are integrating over $\Gamma_6$, we have that $N_1^\star
\sim N_2^\star$. Moreover, we may assume that $N_{1}^{\star}\sim
N_{2}^{\star} \gtrsim N$, otherwise $M_{6}-\prod_{i=1}^{6}m_{i}
\equiv 0$ and \eqref{perturbation} follows trivially.

Since $\frac{1}{m(N_{1}^{\star})(N_{1}^{\star})^{1-}} \lesssim
\frac{1}{N^{1-}}$ we have that
$$\int_{\Gamma_{6}} \prod_{j=1}^{6}\widehat{u_j}(k_{j}) \lesssim
\frac{({N_1^\star})^{0-}}{N^{1-}}\int (\widehat{J
Iu_1^\star})\prod_{j=2}^{6}\widehat{u_j^\star} \lesssim
\frac{({N_1^\star})^{0-}}{N^{1-}}\|Iu_1^\star\|_{H^{1}_x}
\prod_{j=2}^{6}\|u_j^\star\|_{L_{x}^{10}}$$ by reversing
Plancherel and applying H\"older's inequality. Moreover by Sobolev
embedding
$$\|u_j^\star\|_{L_{x}^{10}} \lesssim
\|u_j^\star\|_{H^{2/5}_x}$$ and for $s \geq 2/5$ we have
$$\|u_j^\star\|_{H^{2/5}_x} \lesssim \|Iu_j^\star\|_{H^{1+\frac{2}{5}-s}_x} \lesssim \|Iu_j^\star\|_{H^{1}_x}.$$
Thus \eqref{pieceb} follows.
\end{proof}

\subsection*{Step 3: Decay of ${\mathbf{E^{2}(u)}}$}
In order to estimate the decay of $E^{2}(u)$ in time we need to
bound the right hand side of \eqref{increment}. Towards this aim
we prove a bilinear estimate using an elementary number theory
argument. Let $\eta \in C_{0}^{\infty}$ be a bump function that is
supported on $[-2,2]$ and equals 1 on $[-1,1]$.

\begin{prop}\label{counting}
Let $\phi_1$ and $\phi_2$ be $\lambda-$periodic functions whose
Fourier transforms are supported on $\{k:|k| \sim N_1\}$ and
$\{k:|k| \sim N_2\}$ respectively with $N_1 >> N_2$. Then
$$\|\eta(t) (U_{\lambda}(t)\phi_1)\; \eta(t) (U_{\lambda}(t)\phi_2) \|_{L_{t}^{2}L_{x}^{2}}
\lesssim C(\lambda, N_1)\|\phi_1\|_{L^{2}} \|\phi_2\|_{L^{2}}$$
where
\[C(\lambda ,N_1)= \left\{\begin{array}{ll}
1, & \mbox{if $N_1 \leq 1$}\\
(\frac{1}{N_1}+\frac{1}{\lambda})^{\frac{1}{2}}, & \mbox{if
$N_1>1$}
\end{array}
\right\}.\]
Moreover, 
$$\|\phi_1 \; \phi_2 \|_{L_{t}^{2}L_{x}^{2}}
\lesssim C(\lambda, N_1)\|\phi_1\|_{X^{0, 1/2+}} \|\phi_2\|_{X^{0,1/2+}}.$$

\end{prop}

\subsection*{Remark} We observe that as
$\lambda \rightarrow \infty$ we recover the refined
bilinear estimate in ${\mathbb R}$ with constant 
$C(N_1) = (\frac{1}{N_1})^{1/2}$, see, for example, \cite{bo3, ckstt1, ot1998}.

\begin{proof}
Let $\psi$ be a positive even Schwartz function such that $\psi =
\hat{\eta}$. Then we have
\begin{align}
& B = \|\eta(t) (U_{\lambda}(t)\phi_1)\; \eta(t) (U_{\lambda}(t)\phi_2) \|_{L_{t}^{2}L_{x}^{2}} \nonumber \\
& = \left\|\int_{k=k_1+k_2,\;\tau=\tau_1+\tau_2}
\widehat{\phi_1}(k_1) \widehat{\phi_2}(k_2) \psi(\tau_1-k_{1}^{2})
\; \psi(\tau_2-k_{2}^{2}) \; (dk_1)_{\lambda} \; (dk_2)_{\lambda}
\;
d\tau_{1} \; d\tau_{2} \right\|_{L_{\tau}^{2}L_{k}^{2}} \nonumber \\
& \lesssim \left\| \left(\int_{k=k_1+k_2} \widetilde{\psi}(\tau - k_{1}^{2} - k_2^{2}) \; (dk_1)_{\lambda} \; (dk_2)_{\lambda} \right)^{1/2} \times\right. \;\nonumber \\
& \left.\times \left( \int_{k=k_1+k_2} \widetilde{\psi}(\tau -
k_1^2 - k_2^2) \; |\widehat{\phi}_{1}(k_1)|^{2} \;
|\widehat{\phi}_{2}(k_2)|^{2} \; \; (dk_1)_{\lambda} \;
(dk_2)_{\lambda} \right)^{1/2} \right\|_{L_{\tau}^{2}L_{k}^{2}},
\label{bilCS}
\end{align}
where to obtain \eqref{bilCS} we used Cauchy-Schwartz and the
following definition of $\widetilde{\psi} \in \mathcal{S}$
$$\int_{\tau = \tau_1 + \tau_2}
\psi(\tau_1 - k_1^2) \; \psi(\tau_2 - k_2^2) \; d\tau_1 d\tau_2 =
\widetilde{\psi}(\tau - k_1^2 - k_2^2).$$ An application of
H\"{o}lder gives us the following upper bound on \eqref{bilCS}
\begin{equation} \label{bilHol}
M \left\| \left( \int_{k=k_1+k_2} \widetilde{\psi}(\tau - k_1^2 -
k_2^2) \; |\widehat{\phi}_{1}(k_1)|^{2} \;
|\widehat{\phi}_{2}(k_2)|^{2} \; \; (dk_1)_{\lambda} \;
(dk_2)_{\lambda} \right)^{1/2} \right\|_{L_{\tau}^{2}L_{k}^{2}},
\end{equation}
where
$$ M =\left\| \int_{k=k_1+k_2} \widetilde{\psi}(\tau - k_1^2 - k_2^2) \; (dk_1)_{\lambda}
\; (dk_2)_{\lambda} \right\|_{L_{\tau}^{\infty}L_{k}^{\infty}}
^{1/2}.$$ Now by integration in $\tau$ followed by Fubini in
$k_1$, $k_2$ and two applications of Plancharel we have
$$ \left\| \left( \int_{k=k_1+k_2} \widetilde{\psi}(\tau - k_1^2 - k_2^2) \;
|\widehat{\phi}_{1}(k_1)|^{2} \; |\widehat{\phi}_{2}(k_2)|^{2} \;
\; (dk_1)_{\lambda} \; (dk_2)_{\lambda} \right)^{1/2}
\right\|_{L^{2}_{k, \tau}} \lesssim \|\phi_1\|_{L^2_x}
\|\phi_2\|_{L^2_x},$$ which combined with \eqref{bilCS},
\eqref{bilHol} gives
\begin{equation} \label{bilest}
B \lesssim M \|\phi_1\|_{L^2_x} \|\phi_2\|_{L^2_x}.
\end{equation}

We find an upper bound on $M$ as follows:
\begin{equation} \label{M}
M \lesssim \left( \frac{1}{\lambda} \sup_{\tau, k} \#S
\right)^{1/2},
\end{equation}
where
$$ S = \{ k_1 \in \frac{1}{\lambda}\mathbb{Z} \; | \; |k_1| \sim N_1, \; |k - k_1| \; \sim N_2, \; k^2 - 2k_1(k-k_1) = \tau + O(1) \},$$
and $\#S$ denotes the number of elements of $S$.
When $N_1 \leq 1$ then $\#S \lesssim O(\lambda)$, which implies
$C(\lambda, N_1) \lesssim 1$. When $N_1 > 1$ \, rename $k_1 = z$.
Then
$$ S = \{ z \in \frac{1}{\lambda}\mathbb{Z}\; | \; |z| \sim N_1, \; |k - z| \; \sim N_2, \;
k^2 - 2z(k-z) = \tau + O(1) \}.$$ Let $z_0$ be an element of $S$
i.e.
\begin{equation} \label{z0fr}
|z_0| \sim N_1, \; \; |k-z_0| \sim N_2,
\end{equation}
and
\begin{equation} \label{z0}
k^{2} + 2 z_0^{2} - 2 k z_0 = \tau + O(1).
\end{equation}
In order to obtain an upper bound on $\#S$, we shall count the
number of $\bar{z}$'s $\in \frac{1}{\lambda}{\mathbb{Z}},$ such
that $z_0 + \bar{z} \in S$. Thus
\begin{equation} \label{zmfr}
|z_0 + \bar{z}| \sim N_1, \; \; |z_0 +\bar{z} - k| \sim N_2,
\end{equation}
and
\begin{equation}  \label{zm}
k^{2} + 2 \left(z_0 + \bar{z}\right)^{2} - 2 k \left(z_0 +
\bar{z}\right)= \tau + O(1).
\end{equation}
However by \eqref{z0} we can rewrite the left hand side of
\eqref{zm} as follows
\begin{align*}
k^{2} + 2 \left(z_0 + \bar{z}\right)^{2} - 2 k \left(z_0 +
\bar{z}\right) & =
k^{2} + 2 z_0^{2} + 2 \left(\bar{z}\right)^{2} + 4 z_0\bar{z} - 2 k z_0 - 2 k \bar{z}\\
& =  \tau + O(1) + 2 \bar{z}^{2} + 4 z_0\bar{z} - 2 k \bar{z}.
\end{align*} Hence it suffices to count $\bar{z}$'s $\in \frac{1}{\lambda}\Bbb Z$ satisfying
\eqref{zmfr} and such that
\begin{equation} \label{z0m}
\bar{z}^{2} + 2\bar{z} (z_0 - \frac{k}{2}) = O(1),
\end{equation}
where $z_0$ satisfies \eqref{z0fr} - \eqref{z0}.

By \eqref{z0fr} and \eqref{zmfr} we have that
$$\left|\bar{z}\right| = \left|\left(\bar{z} + z_0 - k \right)- z_0 + k\right| \lesssim N_2 + N_2.$$
Hence
\begin{equation} \label{ballN2}
|\bar{z} | \lesssim N_2 << N_1.
\end{equation}
On the other hand, by \eqref{z0fr} we have that
\begin{equation}\label{sizeN1}|z_0 - \frac{k}{2}| \sim N_1.\end{equation} Now rewrite
\eqref{z0m} as follows
\begin{equation} \label{z0mrev1}
\bar{z} \left( \bar{z} + 2 (z_0 - \frac{k}{2}) \right) = O(1).
\end{equation}
Therefore, using \eqref{ballN2}-\eqref{sizeN1}, the estimate
\eqref{z0mrev1} implies that
$$\left|\bar{z}\right| N_1 = O(1).$$
Since $\bar{z} \in \frac{1}{\lambda}\Bbb Z$ this implies that the
number of $\bar{z}$'s of size $1/N_1$ is $\lambda/N_1$. Hence
$$\#S \lesssim 1 + \frac{\lambda}{N_1},$$
which combined with \eqref{M} gives
$$ M \leq \left(\frac{1}{N_1} + \frac{1}{\lambda}\right)^{1/2},$$
and therefore
$$C(\lambda, N_1) \leq \left(\frac{1}{N_1}+\frac{1}{\lambda}\right)^{\frac{1}{2}}.$$
%In the case that $\lambda \lesssim N_{1}$ then obviously
%$$\#S \lesssim 2$$
%and
%$$C(\lambda, N_1) \leq \left (\frac{2}{\lambda}\right )^{\frac{1}{2}}.$$
\end{proof}

Now we are ready to prove desired decay of $E^{2}(u)$ which
follows from the following proposition:

\begin{prop}\label{edec}
For any $\lambda$-periodic Schwartz function u, with period
$\lambda \geq N$, and any $\delta$ given by the local theory, we have that
$$|\int_{0}^{\delta}\Lambda_{10}(M_{10};u(t))| \lesssim 
\lambda^{0+} N^{-5/2+}\|Iu\|_{X^{1,1/2+}}^{10},$$ for $s > 11/28.$
\end{prop}

\begin{proof}
We perform a dyadic decomposition as in
Proposition \ref{edif}, and we borrow the same notation. 
From now on we do not keep track of all the different $\lambda^{k \epsilon}$, 
$k \in \Bbb R$, and we just write $\lambda^{0+}$

We observe that $M_{10}$ is bounded as an elongation of the
bounded multiplier $M_6$. Since the multiplier $M_{10}$ vanishes on $\Gamma_{10}$ when all frequencies 
are smaller than $N$, we can assume that there are $N_1^{*}, N_2^{*} \gtrsim N$.
Now we divide the proof in two cases. 
\\
\\
{\bf Case 1}. $N_{1}^{\star} \sim N_{2}^{\star} \sim N_{3}^{\star}
\gtrsim N$
\\
\\
Since
$\frac{1}{m(N_{1}^{\star})m(N_{2}^{\star})m(N_{3}^{\star})(N_{1}^{\star}N_{2}^{\star}N_{3}^{\star})^{1-}}
\lesssim \frac{1}{N^{3-}}$ we have
\begin{align}
|\int_{0}^{\delta}\int M_{10}\prod_{j=1}^{10}\hat{u}_{j}|
&\lesssim \frac{(N_1^\star)^{0-}}{N^{3-}}
\|JIu_{1}^\star\|_{L_{t}^{6}L_{x}^{6}}\|JIu_{2}^\star\|_{L_{t}^{6}L_{x}^{6}}\|JIu_{3}^\star\|_{L_{t}^{6}L_{x}^{6}}
\prod_{j=4}^{10}\|u_{j}^\star\|_{L_{t}^{14}L_{x}^{14}} \nonumber
\\ & \lesssim
\lambda^{0+}\frac{(N_1^\star)^{0-}}{N^{3-}}\|Iu\|_{X^{1,1/2+}}^{3}\|u\|_{X^{(\frac{1}{2}-\frac{3}{14})+,1/2+}}^{7}
\label{useinter}\\ & \lesssim
\lambda^{0+}(N_1^\star)^{0-}N^{-3+}\|Iu\|_{X^{1,1/2+}}^{10}\label{spaces}
\end{align} where in order to obtain \eqref{useinter} we use
\eqref{int1} and to obtain \eqref{spaces} we use the fact that for
$s>2/7$ the following inequality holds
$$\|u\|_{X^{(\frac{1}{2}-\frac{3}{14})+,1/2+}}=\|u\|_{X^{2/7+,1/2+}}\lesssim \|Iu\|_{X^{1,1/2+}}.$$
\\
{\bf Case 2}. $N_{1}^{\star} \sim N_{2}^{\star} \gg N_{3}^{\star}$
\\
\\
Since $\delta<1$ we can insert $\eta(t)$ where $\eta$ is
a bump function supported in $[-1/2,3/2]$ and equals 1 in $[0,1]$.
Also notice that
$\frac{1}{m(N_{1}^{\star})m(N_{2}^{\star})(N_{1}^{\star}N_{2}^{\star})^{1-}}
\lesssim \frac{1}{N^{2-}}$. Hence using the Cauchy-Schwartz
inequality we get:
\begin{align}
|\int_{0}^{\delta}\eta(t)\int M_{10}\prod_{j=1}^{10}\hat{u}_{j}| 
& \lesssim \frac{(N_1^\star)^{0-}}{N^{2-}} 
\|\eta JIu_1^\star u_3^\star\|_{L^2_tL^2_x}
\|JIu_2^\star \prod_{4}^{10}u_i^\star \|_{L^2_tL^2_x}\nonumber \\
& \lesssim \frac{(N_1^\star)^{0-}}{N^{2-}} 
(\frac{1}{\lambda} + \frac{1}{N_{1}^{*}})^{\frac{1}{2}} 
\|JIu_1^\star\|_{X^{0,1/2+}}\| u_3^\star\|_{X^{0,1/2+}} 
\|JIu_2^\star \prod_{4}^{10}u_i^\star \|_{L^2_tL^2_x}\label{usebil}\\
& \lesssim \frac{(N_1^\star)^{0-}}{N^{5/2-}}\|Iu_1^\star
\|_{X^{1,1/2+}}\|u_3^\star\|_{X^{0,1/2+}}\|JIu_2^\star\|_{L^4_tL^4_x}\prod_{j=4}^{10}
\|u_j^\star\|_{L^{28}_tL^{28}_x}\label{useholder}\\
& \lesssim \lambda^{0+}\frac{(N_1^\star)^{0-}}{N^{5/2-}}\|Iu
\|_{X^{1,1/2+}}^{3}\prod_{j=4}^{10} \|u_j\|_{X^{1/2-3/28,1/2+}}\label{useinter2}\\
& \lesssim \lambda^{0+}(N_1^\star)^{0-}N^{-5/2+}\|Iu\|_{X^{1,1/2+}}^{10}
\nonumber
\end{align}
where to obtain \eqref{usebil} we use Proposition \ref{counting},
to obtain \eqref{useholder} we use H\"older's inequality and the
facts that $\lambda \geq N$ and $N_1^* \gtrsim N$, to obtain \eqref{useinter2} we use
\eqref{int1}, and in the last line we use that $s>11/28$.
\end{proof}

Now we are ready to prove our main theorem on ${\Bbb T}^{1}$.
\subsection*{Step 4: Proof of Theorem 1}
\begin{proof}Let $u_0 \in H^s$, where $4/9 < s < 1/2$. By definition,
$\frac{|m(k)|}{| k|^{s-1}} \lesssim N^{1-s}$, hence
$$\| \partial_{x} Iu_{0}^{\lambda}\|_{2}=\|m(k)| k| \widehat{u_{0}^{\lambda}}\|_{2}=
\| \frac{|m(k)|}{| k|^{s-1}}|
k|^{s}\widehat{u_{0}^{\lambda}}\|_{2} \leq N^{1-s} \|| k|^{s}
\widehat{u_{0}^{\lambda}}\|_{2} \lesssim
\frac{N^{1-s}}{\lambda^{s}} \|u_{0}\|_{\dot{H}^{s}}.$$ By the
Gagliardo-Nirenberg inequality we have:
$$\| Iu_{0}^{\lambda}\|_{6}^{6} \lesssim \| \partial_{x} Iu_{0}^{\lambda}\|_{2}^{2}\|Iu_{0}^{\lambda}\|_{2}^{4} \lesssim
 \frac{N^{2-2s}}{\lambda^{2s}} \|u_{0}\|_{\dot{H}^{s}}^{2} \|u_{0}\|_{2}^{4}.$$
Combining the two inequalities above with the definition of the
first modified energy we obtain
$$E^1(u_{0}^{\lambda})=
\frac{1}{2}\int_{\Bbb
R}|\partial_{x}Iu_{0}^{\lambda}|^{2}dx+\frac{1}{6} \int_{\Bbb
R}|Iu_{0}^{\lambda}|^{6}dx \lesssim \frac{N^{2-2s}}{\lambda^{2s}}
\|u_{0}\|_{\dot{H}^{s}}^{2}.$$ We choose $\lambda \sim N^{
\frac{1-s}{s}}$, which for $s < 1/2$ implies $\lambda \gtrsim N$.
Then we have that $E^1(u_{0}^{\lambda}) \lesssim 1$. Therefore
\begin{equation} \label{init}
\|Iu_{0}^{\lambda}\|_{H^{1}}^{2}\lesssim 1,
\end{equation}
which allows us to apply Propositions \ref{lwp}, \ref{edif}. The
main idea of the proof is that if in each step of the iteration we
have the same bound for $\|Iu_{0}^{\lambda}\|_{H^{1}}$ then we can
iterate again with the same timestep. Because we iterate the 
$u^{\lambda}$ solutions, the interval of local existence
is of order $\frac{1}{\lambda^{a\epsilon}}$. But this creates 
no additional problem since $\lambda$ is given in terms
of $N$. Moreover $N$ is a fixed number and thus the interval of 
local existence does not shrink. This is the reason
why we don't write down the $\lambda^{0+}$ dependence 
in the following string of inequalities. Now, by Lemma \ref{fundcal},
Proposition \ref{edif}, Proposition \ref{edec} and \eqref{init} we
have that
$$E^{2}(u^{\lambda}(\delta)) \lesssim E^{2}(u^{\lambda}(0))+CN^{-5/2+} \lesssim
E^{1}(u^{\lambda}(0))+\frac{1}{N}\|Iu_{0}^{\lambda}\|_{H^{1}}^{6}+CN^{-5/2+}
\lesssim 1,$$ 
for $N$ large enough. Combining the inequality above
with  Proposition \ref{edif} we have,
$$\| \partial_{x}Iu(\delta)^{\lambda}\|_{2}^{2} \lesssim 
E^{2}(u^{\lambda}(\delta)) + O(\frac{1}{N}) \lesssim 1.$$
Thus
$$\|Iu^{\lambda}(\delta)\|_{H^{1}}\lesssim 1$$
and we can continue the solution in $[0,M\delta]=[0,T]$ as long as
$T\ll N^{5/2-}$. Hence
$$\|Iu^{\lambda}(T)\|_{H^{1}} \lesssim 1$$
for all $T\ll N^{5/2-}$. From the definition of $I$ this implies
that
$$\|u^{\lambda}(T)\|_{H^{s}} \lesssim 1$$
for all $T\ll N^{5/2-}$. Undoing the scaling we have that
$$\|u(T)\|_{H^{s}} \lesssim C_{N,\lambda}$$
for all $T\ll \frac{ N^{5/2-}}{\lambda^{2}}$. But $\frac{
N^{5/2-}}{\lambda^{2}} \sim N^{\frac{9s-4}{2s}}$ goes to infinity
as $N\rightarrow \infty$ since $s>4/9.$
\end{proof}

\section{The I-method and the proof of Theorem \ref{main2D}.}
%Consider the semilinear Schr\"odinger initial value problem
%\begin{align} \label{ivp2}
%&iu_{t}+\Delta u-|u|^{2}u=0\\ \label{bc} &u(x,0)=u_{0}(x)\in
%H^{s}({\Bbb T^2 }),t\in {\Bbb R},
%\end{align} for $0 < s <1.$ We recall that a solution $u$ to
%(\ref{ivp2})-(\ref{bc}), satisfies the following two conservation
%laws:

%We will work with $\lambda-$periodic functions, thus we define
%$(dk)_{\lambda}$ to be the normalized counting measure on
%$\frac{1}{\lambda}\Bbb Z^2$:
%$$\int a(x)(dk)_{\lambda}=\frac{1}{\lambda^2}\sum_{k \in \frac{1}{\lambda}\Bbb Z}a(k).$$
%Then define the Fourier transform of $f(x) \in L_{x\in
%[0,\lambda]}^{1}$ by
%$$\hat{f}(k)=\int_{0}^{\lambda}e^{-2\pi ikx}f(x)dx$$
%For an appropriate class of functions we know by the Fourier
%inversion formula that:
%$$f(x)=\int e^{2\pi ikx}\hat{f}(k)(dk)_{\lambda}.$$

In this section we present the proof of Theorem \ref{main2D}. We
will focus on the analysis of the $\lambda$-periodic problem:
\begin{align}
&iu_{t}+ \Delta u -|u|^{2}u=0 \label{lambdaivp}\\
&u(x,0)=u_{0}(x)\in H^{s}({\Bbb T_\lambda^{2} }),t\in {\Bbb R}
\label{lambdabc},
\end{align}
for $0<s<1$. As for the one dimensional case, we will use the
$I$-method. The definitions of the multiplier $m$ and of the
operator $I$ are the same as in the 1D context. Precisely, given a
parameter $N>>1$ to be chosen later, we define $m(k)$ to be the
following multiplier:

\[m(k)= \left\{\begin{array}{ll}
1 & \mbox{if $|k|<N$}\\
\left(\frac{|k|}{N}\right)^{s-1} & \mbox{if $|k|>2N$}
\end{array}
\right.\] Then $I:H^{s}\rightarrow H^{1}$  will be defined as the
following multiplier operator:
$$(\widehat{Iu})(k)=m(k)\hat{u}(k).$$ This operator is
smoothing of order $1-s$, that is:
\begin{equation}\label{Iprop}
\|u\|_{s_{0},b_{0}}\lesssim \|Iu\|_{s_{0}+1-s,b_{0}}\lesssim
N^{1-s}\|u\|_{s_{0},b_{0}}
\end{equation}
for any $s_{0},b_{0}\in{\Bbb R}$.

The first modified energy of a function $u \in H^s$ is defined by
$$E^1(u)=E(Iu),$$ where we recall that
$$E(u)(t)=\frac{1}{2}\int
|\nabla u(t)|^{2}dx+\frac{1}{4}\int |u(t)|^{4}dx.$$

Our main objective will be to prove that the modified energy of a
solution to the $\lambda$-periodic problem
\eqref{lambdaivp}-\eqref{lambdabc} is almost conserved in time.

%Indeed one can easily prove the following bound:
%$$\|u(\cdot,t)\|_{H^{s}(\Bbb T^2)} \lesssim E^1(u)(t) + \|u_0\|_{L^2(\Bbb T^2)},$$
%for all solutions $u$ to (\ref{ivp2})-(\ref{bc}). Hence, it is
%enough to control $E^1(u)(t)$, in order to control the growth of
%the $H^s$ norm of $u$.

%Now, we collect some known estimates, which will be often used in
%our arguments.
%First, we recall the following  Strichartz type estimate (see
%\cite{bo}), which holds for $\lambda$-periodic functions, whose
%Fourier transform is supported on $\{k \in
%\frac{1}{\lambda}\mathbb{Z}^2: |k| \sim N\}$ :
%\begin{equation} \label{L44}
%\|u\|_{L_{t}^{4}L_{x}^{4}} \lesssim \lambda^{0+} N^{0+}
%\|u\|_{X^{0,1/2+}}
%\end{equation}
%We also have the following ``energy type'' estimate:
%\begin{equation} \label{energy}
%\|u\|_{L_{t}^{\infty}L_{x}^{\infty}} \lesssim \lambda^{0+}
%\|u\|_{X^{1/2+,1/2+}}.
%\end{equation}
%Hence, by interpolation, we get
%\begin{equation} \label{interpol}
%\|u\|_{L_{t}^{p}L_{x}^{p}} \lesssim \lambda^{0+}N^{0+}
%\|u\|_{X^{\alpha(p),1/2+}},
%\end{equation}
%with $\alpha(p)=\frac{1}{2}(1-4/p)$, and $4 \leq p \leq \infty.$

\

We now proceed to present the steps leading to the proof of
Theorem \ref{main2D}.

\subsection*{Step 1: Local well-posedness for the ${\mathbf{I}}$-system}
The first step towards the proof of Theorem \ref{main2D} is to
apply the $I$-operator to \eqref{ivp1} and prove a local
well-posedness result for the $I$-initial value problem
\begin{align}\begin{split} \label{2Divp}
&iIu_{t}+I\Delta u-I(|u|^{2}u)=0\\
&Iu(x,0)=Iu_{0}(x)\in H^{s}({\Bbb T_{\lambda}^{2} }),t\in {\Bbb R}.
\end{split}
\end{align}
This is the content of the next proposition.

\begin{prop}\label{2Dlwp}
Consider the $I$-initial value problem \eqref{2Divp}, $\lambda$=1. If
$\|Iu_{0}\|_{H^{1}} \lesssim 1$, then \eqref{2Divp} is locally
well-posed for any $s>\frac{7}{20}$ in $[0,\delta]\sim [0,1]$.
\end{prop}
\begin{proof}
By Duhamel's formula as in Proposition \ref{lwp} we have that
$$\|Iu\|_{X^{1,1/2+}} \lesssim 
\|Iu_{0}\|_{H^{1}}+\delta^{\frac{1}{20}-\epsilon}\|I(|u|^{2}u)\|_{X^{1,-\frac{9}{20}}}$$
By the ``invariant lemma'' in \cite{ckstt2} we know that the estimate
$$\|I(|u|^{2}u)\|_{X^{1,-\frac{9}{20}}}\lesssim \|Iu\|_{X^{1,1/2+}}^{3}$$
is implied by 
\begin{equation} \label{tril}
\||u|^{2}u\|_{X^{s,-\frac{9}{20}}}\lesssim \|u\|_{X^{s,1/2+}}^{3}.
\end{equation}
Since the rest of the proof is identical to Proposition \ref{lwp}
we briefly recall how \eqref{tril} can be obtained. 
Note that this estimate was first proved by Bourgain \cite{bo1} 
and his proof is based upon a local variant of the well known 
periodic Strichartz estimate
\begin{equation} \label{2DL4loc}
\|u\|_{L_{t}^{4}L_{x}^{4}} \lesssim 
N^{s_{1}}\left( \sum_{k \in Q}\int d\tau(1+|\tau-|k|^{2}|)^{2b_{1}}
|\hat{u}(k,\tau)|^{2} \right)^{\frac{1}{2}}
\end{equation}
where $b_{1}>\frac{1-\min(1/2,s_{1})}{2}$.
Here we omit the proof and refer the reader to \cite{bo1} for details. 
\end{proof}

%%%%%%%%%%%%%%%%%%%%%%%%%%%%%%%%%%%%%%%%%%%%%%%%%%%%%%%%%%%%%%%%%%%%
{} %end of the comment

%%%%%%%%%%%%%%%%%%%%%%%%%%%%%%%%%%%%%%%%%%%%%%%%%%%%%%%%%%%%%%%%%%%%%%%%%%%%%%%%%%

\begin{cor}\label{n2Dlwp}
Consider the $I$-initial value problem \eqref{2Divp}. If
$\|Iu_{0}\|_{H^{1}} \lesssim 1$, then \eqref{2Divp} is locally
well-posed for any $s>\frac{7}{20}$ in $[0,\delta]\sim [0,\frac{1}{\lambda^{a\epsilon}}]$, where $a$ is a fixed positive integer
($a \sim 100$ will do it).
\end{cor}
\begin{proof}
The proof is identical to the argument of Bourgain in \cite{bo1}, 
the only difference being the fact that \eqref{2DL4loc} holds true
for the $\lambda$-periodic solutions, but now with a factor of 
order $\lambda^{0+}$ on the right hand side. 
Again the fixed power of $\delta$ in the local theory, enables us
 to prove the local well-posedness result for 
the $\lambda$-periodic problem in the interval 
$[0,\frac{1}{\lambda^{a\epsilon}}]$.
\end{proof}
%As an immediate consequence of (\ref{L44}), we get the following
%bilinear estimate, for $\lambda$-periodic functions $v_1$ and
%$v_2$, supported on frequency space around $N_1$ and $N_2$
%respectively:
%\begin{equation} \label{L22}
%\|v_1 v_2\|_{L_{t}^{2}L_{x}^{2}} \lesssim \lambda^{0+}
%(N_1N_2)^{0+} \|v_1\|_{X^{0,1/2+}}\|v_2\|_{X^{0,1/2+}}.
%\end{equation}

\subsection*{Step 2: Decay of the first energy} In order to prove
that the first energy is almost conserved, we will use the
bilinear estimate for $\lambda$-periodic functions stated in
Proposition \ref{2Dcounting}. Its proof is based on some
number theoretic facts that we recall in the following
three lemmas; see also related  estimates in the work of Bourgain \cite{boirr}.

The following lemma is known as {\bf Pick's Lemma} \cite{pick}:
\begin{lem}\label{Pick}
Let $Ar$ be the area of a simply connected lattice polygon.
Let $E$ denote the number of lattice points on the polygon edges and
$I$ the number of lattice points in the interior of the polygon. Then
$$ Ar = I + \frac{1}{2} E - 1.$$
\end{lem}

\begin{lem}\label{latticepts}
Let $\mathcal{C}$ be a circle of radius $R$. If $\gamma$ is an arc
on $\mathcal{C}$ of length $|\gamma| <
\left(\frac{3}{4}R\right)^{1/3}$, then $\gamma$ contains at most 2
lattice points.
\end{lem}

\begin{proof}
We prove the lemma by contradiction. Assume that there are 3 lattice
points $P_1$, $P_2$ and $P_3$ on an arc $\gamma = AB$ of $\mathcal{C}$, 
and denote by $T(P_1,P_2,P_3)$ the triangle with vertices $P_1$, $P_2$ and $P_3$. 
Then, by Lemma
\ref{Pick} we have
$$\mbox{Area of}\;T(P_1,P_2,P_3) = I + \frac{1}{2}E - 1 \geq I + \frac{3}{2} - 1 = I + \frac{1}{2} \geq \frac{1}{2}.$$
We shall prove that under the assumption that $|\gamma| <
\left(\frac{3}{4}R\right)^{1/3}$, then
\begin{equation} \label{area}
\mbox{Area of}\; T(P_1,P_2,P_3) < \frac{1}{2},
\end{equation}
hence $\gamma$ must contain at most two lattice points.

We observe that (see Figure \ref{fig:triangle})
\begin{figure}
\centering
\scalebox{0.5}{
        \epsfig{file=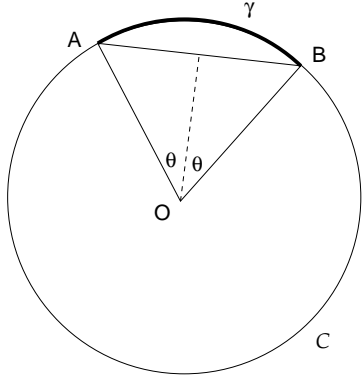}
}
\caption{Triangle area.
}
\label{fig:triangle}
\end{figure}

$$ \mbox{Area of the sector} \; ABO = R^2 \theta,$$
$$ \mbox{Area of the triangle} \; ABO = R^2 \sin \theta \cos \theta.$$
Hence, for any $P_1$, $P_2$, $P_3$ on $\gamma$ we have
\begin{equation} \label{areacal}
\mbox{Area of} \; T(P_1,P_2,P_3) \leq R^2 \theta -  R^2 \sin
\theta \cos \theta = R^2 (\theta - \frac{1}{2}\sin (2\theta)).
\end{equation}
One can easily check that
\begin{equation} \label{afterTay}
\theta - \frac{1}{2}\sin (2\theta) \leq   \frac{2}{3} \theta^3.
\end{equation}
Thus \eqref{areacal}, \eqref{afterTay} and the fact that $|\gamma|
= R \theta$ imply that
$$ \mbox{Area of} \; T(P_1,P_2,P_3) \leq \frac{2}{3} R^{2} \theta^3
= \frac{2}{3} R^2 (|\gamma| R^{-1})^{3} < \frac{1}{2},$$ where to
obtain the last inequality we used the assumption that $|\gamma| <
\left(\frac{3}{4}R\right)^{1/3}$. Therefore \eqref{area} is
proved.

\end{proof}

Also we recall the following result of Gauss, see, for example \cite{ai} 
\begin{lem}\label{Gauss}
Let $K$ be a convex domain in ${\mathbb{R}}^2$. If
$$ N(\lambda) = \#\{ {\mathbb{Z}}^2 \cap \lambda K \},$$
then, for $\lambda >> 1$
$$ N(\lambda) = \lambda^2 |K| + O(\lambda),$$
where $|K|$ denotes the area of $K$ and $\# S$ denotes the number of
points of a set $S$.
\end{lem}

Now we are ready to present the bilinear estimate for $\lambda$-periodic functions.
Let $\eta \in C_{0}^{\infty}$ be a bump function that is supported
on $[-2,2]$ and equals 1 on $[-1,1]$.

\begin{prop}\label{2Dcounting}
Let $\phi_1$ and $\phi_2$ be $\lambda-$periodic functions,
$\lambda \geq 1$, whose Fourier transforms are supported on
$\{k:|k| \sim N_1\}$ and $\{k:|k| \sim N_2\}$ respectively, with
$N_1 >> N_2 > 1$.
\begin{enumerate}
\item[(a)] Then
\begin{equation}\label{bila}
\|\eta(t) (U_{\lambda}(t)\phi_1)\;\eta(t) (U_{\lambda}(t)\phi_2)
\|_{L_{t}^{2}L_{x}^{2}} \lesssim (\lambda N_2)^{\epsilon}
\|\phi_1\|_{L^{2}} \|\phi_2\|_{L^{2}},\end{equation} for any
$\epsilon > 0$.
Hence 
$$\|\phi_1\;\phi_2\|_{L_{t}^{2}L_{x}^{2}} 
\lesssim (\lambda N_2)^{\epsilon}
\|\phi_1\|_{X^{0,1/2+}} \|\phi_2\|_{X^{0,1/2+}}.$$

\item[(b)] Moreover, if $\lambda >> 1$ then
\begin{equation}\label{bilb}\|\eta(t) (U_{\lambda}(t)\phi_1)\; \eta(t) (U_{\lambda}(t)\phi_2)
\|_{L_{t}^{2}L_{x}^{2}} \lesssim \left( \frac{1}{\lambda} +
\frac{N_2}{N_1}\right)^{1/2} \|\phi_1\|_{L^{2}}
\|\phi_2\|_{L^{2}}.\end{equation}
\end{enumerate}
\end{prop}

\subsection*{Remark} Note that as $\lambda \rightarrow \infty$ we recover 
the improved bilinear Strichartz with constant $C(N_1, N_2) = (\frac{N_2}{N_1})^{1/2}$. 

\begin{proof}
(a) Let $\psi$ be a positive even Schwartz function such that $\psi = \hat{\eta}$.
Then we have
\begin{align}
B & = \|\eta(t) (U_{\lambda}(t)\phi_1)\;\eta(t) (U_{\lambda}(t)\phi_2) \|_{L_{t}^{2}L_{x}^{2}} \nonumber \\
& = \|\int_{k=k_1+k_2,\;\tau=\tau_1+\tau_2} \widehat{\phi_1}(k_1) \widehat{\phi_2}(k_2)
\psi(\tau_1-k_{1}^{2}) \; \psi(\tau_2-k_{2}^{2}) \; (dk_1)_{\lambda} \; (dk_2)_{\lambda} \;
d\tau_{1} \; d\tau_{2} \|_{L_{\tau}^{2}L_{k}^{2}} \nonumber \\
& \lesssim \| \left(\int_{k=k_1+k_2} \widetilde{\psi}(\tau -
k_{1}^{2} - k_2^{2}) \; (dk_1)_{\lambda} \; (dk_2)_{\lambda}
\right)^{1/2}
\times \nonumber \\
& \times \left( \int_{k=k_1+k_2} \widetilde{\psi}(\tau - k_1^2 -
k_2^2) \; |\widehat{\phi}_{1}(k_1)|^{2} \;
|\widehat{\phi}_{2}(k_2)|^{2} \; \; (dk_1)_{\lambda} \;
(dk_2)_{\lambda} \right)^{1/2} \|_{L_{\tau}^{2}L_{k}^{2}}
\label{2DbilCS}
\end{align}
where to obtain \eqref{2DbilCS} we used Cauchy-Schwartz and the
following definition of $\widetilde{\psi} \in \mathcal{S}$
$$\int_{\tau = \tau_1 + \tau_2}
 \psi(\tau_1 - k_1^2) \; \psi(\tau_2 - k_2^2) \; d\tau_1 \; d\tau_2= \widetilde{\psi}(\tau - k_1^2 - k_2^2).$$
An application of H\"{o}lder gives us the following upper bound on
\eqref{2DbilCS}
\begin{equation} \label{2DbilHol}
M \| \left( \int_{k=k_1+k_2} \widetilde{\psi}(\tau - k_1^2 -
k_2^2) \; |\widehat{\phi}_{1}(k_1)|^{2} \;
|\widehat{\phi}_{2}(k_2)|^{2} \; \; (dk_1)_{\lambda} \;
(dk_2)_{\lambda} \right)^{1/2} \|_{L^{2}_{\tau} L^2_k},
\end{equation}
where
$$ M = \| \int_{k=k_1+k_2} \widetilde{\psi}(\tau - k_1^2 - k_2^2) \; (dk_1)_{\lambda} \; (dk_2)_{\lambda} \|_{L^{\infty}_{\tau} L^{\infty}_k}^{1/2}.$$
Now by integration in $\tau$ followed by Fubini in $k_1$, $k_2$ and two applications of
Plancharel we have
$$ \| \left( \int_{k=k_1+k_2} \widetilde{\psi}(\tau - k_1^2 - k_2^2) \; |\widehat{\phi}_{1}(k_1)|^{2} \;
|\widehat{\phi}_{2}(k_2)|^{2} \; \; (dk_1)_{\lambda} \;
(dk_2)_{\lambda} \right)^{1/2} \|_{L^{2}_{\tau} L^2_k} \lesssim
\|\phi_1\|_{L^2_x} \|\phi_2\|_{L^2_x},$$ which combined with
\eqref{2DbilCS} and \eqref{2DbilHol} gives
\begin{equation} \label{2Dbilest}
B \lesssim M  \|\phi_1\|_{L^2_x} \|\phi_2\|_{L^2_x}.
\end{equation}
We find an upper bound on $M$ as follows:
\begin{equation} \label{2DM}
M \lesssim \left( \frac{1}{{\lambda}^{2}} \sup_{\tau, k} \# S \right)^{\frac{1}{2}},
\end{equation}
where
$$ S = \{ k_1 \in {\mathbb{Z}}^{2}/\lambda \; | \; |k_1| \sim N_1, \; |k - k_1| \; \sim N_2, \;
|k|^2 - 2k_1 \cdot (k-k_1) = \tau + O(1) \},$$ and $\# A$ denotes
the number of lattice points of a set $A$.

For notational purposes, let us rename $k_1 = z$, that is
$$ S = \{ z \in {\mathbb{Z}}^{2}/\lambda \; | \; |z| \sim N_1, \; |k - z| \; \sim N_2, \;
|k|^{2} + 2 |z|^{2} - 2 k \cdot z = \tau + O(1) \}.$$
Let $z_0$ be an element of $S$ i.e.
\begin{equation} \label{2Dz0fr}
|z_0| \sim N_1, \; \; |k-z_0| \sim N_2,
\end{equation}
and
\begin{equation} \label{2Dz0}
|k|^{2} + 2 |z_0|^{2} - 2 k \cdot z_0 = \tau + O(1).
\end{equation}
In order to obtain an upper bound on $\#S$, we shall count the
number of $l$'s $\in {\mathbb{Z}}^{2}$ such that $z_0 +
\frac{l}{\lambda} \in S$ where $z_0$ satisfies \eqref{2Dz0fr} -
\eqref{2Dz0}. Thus such $l$'s must satisfy
\begin{equation} \label{2Dzmfr}
\left|z_0 + \frac{l}{\lambda}\right| \sim N_1, \; \; \left|z_0 +
\frac{l}{\lambda} - k\right| \sim N_2,
\end{equation}
and
\begin{equation}  \label{2Dzm}
|k|^{2} + 2 \left|z_0 + \frac{l}{\lambda}\right|^{2} - 2 k \cdot
(z_0 + \frac{l}{\lambda})= \tau + O(1).
\end{equation}
However by \eqref{2Dz0} we can rewrite the left hand side of \eqref{2Dzm} as follows
\begin{align*}
|k|^{2} + 2 \left|z_0 + \frac{l}{\lambda}\right|^{2} - 2 k \cdot
(z_0 + \frac{l}{\lambda}) & =
|k|^{2} + 2 |z_0|^{2} + 2 \left|\frac{l}{\lambda}\right|^{2} + 4 z_0 \cdot \frac{l}{\lambda} - 2 k \cdot z_0 - 2 k \cdot \frac{l}{\lambda}\\
& =  \tau + O(1) + 2 \left|\frac{l}{\lambda}\right|^{2} + 4 z_0
\cdot \frac{l}{\lambda} - 2 k \cdot \frac{l}{\lambda}.
\end{align*}
Therefore \eqref{2Dzm} holds if
\begin{equation} \label{2Dz0l}
\left|\frac{l}{\lambda}\right|^{2} + 2 \frac{l}{\lambda} \cdot
(z_0 - \frac{k}{2}) = O(1).
\end{equation}
Moreover, \eqref{2Dz0fr} and \eqref{2Dzmfr} yield
$$\left|\frac{l}{\lambda}\right| = \left|\frac{l}{\lambda} + z_0 - k - z_0 + k\right| \lesssim N_2 + N_2,$$
that is
\begin{equation} \label{2DballN2}
|l| \lesssim \lambda N_2.
\end{equation}
Finally we observe that \eqref{2Dz0fr} together with the
assumption that $N_1
>> N_2$ implies that
$$
N_1 \sim N_1 - N_2 \sim \left|\left|\frac{z_0}{2} -
\frac{k}{2}\right| - \left|\frac{z_0}{2}\right|\right| \leq
\left|z_0 - \frac{k}{2}\right| \leq \left|\frac{z_0}{2} -
\frac{k}{2}\right| + \left|\frac{z_0}{2}\right| \sim N_2 + N_1
\lesssim N_1,$$ i.e.
\begin{equation} \label{vector}
\left|z_0 - \frac{k}{2}\right| \sim N_1.
\end{equation}
Hence, it suffices to count the $l's \in {\mathbb{Z}}^{2}$
satisfying \eqref{2Dz0l} and \eqref{2DballN2} where $z_0$ is such
that \eqref{vector} holds.

Let $w = (a,b)$ denote the vector $z_0 - \frac{k}{2}$. Thus we
need to count the number of points in the set $A^{\lambda}$
\begin{equation} \label{Alam}
A^{\lambda} = \{ l \in {\mathbb{Z}}^2 : \; \left|\;|l|^{2} + 2
\lambda l \cdot w\;\right| = O(\lambda^2), \; |l| \lesssim \lambda
N_2, \; |w| \sim N_1 \},
\end{equation}
or equivalently, \begin{equation}\label{B1}A^{\lambda} = \{ (x,y)
\in {\mathbb{Z}}^2 : \; \left|x^2 + y^2 + 2\lambda (ax +
by)\right| \leq c \lambda^2,\; x^2 + y^2 \leq(k_2 \lambda
N_2)^{2},\; a^2 + b^2 \sim N_1^2 \},\end{equation} for some
$c,k_2>0$. Let ${\mathcal{C}}_-^{\lambda}$,
${\mathcal{C}}_{+}^{\lambda}$ be the following circles,
\begin{align*} \begin{split} \label{circ}
& {\mathcal{C}}_-^{\lambda}:\; (x + \lambda a)^2 + (y + \lambda b)^2 = -c\lambda^2 + \lambda^2 (a^2 + b^2)\\
& {\mathcal{C}}_{+}^{\lambda}:\; (x + \lambda a)^2 + (y + \lambda b)^2 = c\lambda^2 + \lambda^2 (a^2 + b^2)\\
\end{split}
\end{align*}
and for any integer $n$, let ${\mathcal{C}}_n^{\lambda}$ be the
circle
$${\mathcal{C}}_n^{\lambda}:\; (x + \lambda a)^2 + (y + \lambda
b)^2 = n + \lambda^2 (a^2 + b^2).$$ Finally, let
${\mathcal{D}}^{\lambda}$ denote the disk
$$ {\mathcal{D}}^{\lambda}:\;x^2 + y^2 \leq (k_2\lambda N_2)^{2}.$$

% Fig lam
\begin{figure}[htbp]
\centering
\includegraphics*[scale=0.7,trim=0in 5.8in 3.5in 0in]{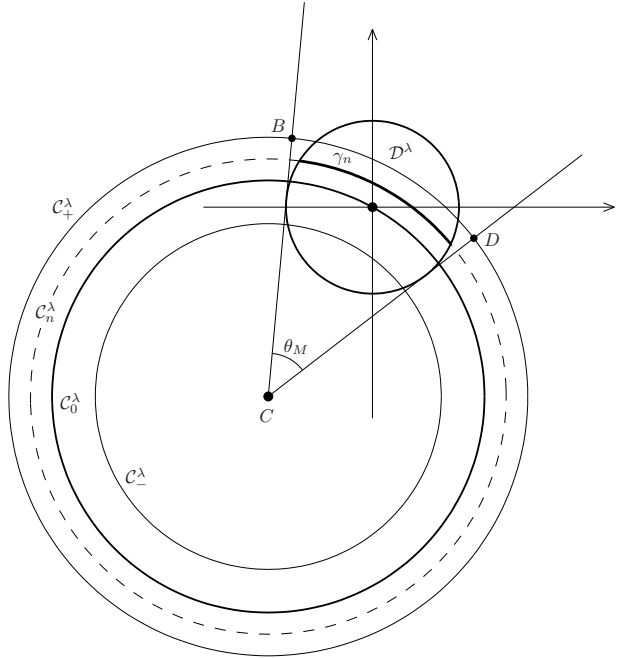}
\caption{Circular sector}
\label{fig:4}
\end{figure}

We need to count the number of lattice points inside
${\mathcal{D}}^{\lambda}$ that are on arcs of circles
${\mathcal{C}}_n^{\lambda}$, with $$ - c \lambda^2 \leq n \leq c
\lambda^2.$$ Precisely, the total number of lattice point in
$A^{\lambda}$ can be bounded from above by
%$$c\lambda^2 \times \mbox{ number of lattice points on } {\mathcal{C}}_{n}^{\lambda}
%\cap {\mathcal{D}}^{\lambda}.$$
%which is well approximated by
\begin{equation} \label{latpts}
2c\lambda^2 \times \#(\mathcal{C}_{n}^{\lambda} \cap
{\mathcal{D}}^{\lambda}).
\end{equation} 

Denote by $\gamma_n$ the arc of circle
$\mathcal{C}_n^{\lambda}$ which is contained in
$\mathcal{D}^{\lambda}$. Notice that (see Figure \ref{fig:4})
\begin{equation}\label{gamman}|\gamma_n| \leq R_M
\theta_{M}\end{equation} where $R_M = \sqrt{c\lambda^2+ k_1
\lambda^2 N_1^2}$ for some constant $k_1>0$, and $\theta_{M}$ is
the angle between the line segment $CB$ and $CD$, which lie along
the tangent lines from $C=(-\lambda a,-\lambda b)$ to the circle
$x^2 + y^2 = (k_2\lambda N_2)^{2}$. Hence,
$$\sin{\theta_{M}} \leq k \frac{N_2}{N_1},$$ for some constant $k>0$. Since $N_1 >> N_2$,
we can assume that $\sin \theta_M > \frac{1}{2}\theta_M$. Hence,
\begin{equation}\label{thetaM}\theta_M < 2k \frac{N_2}{N_1}.\end{equation} In order to count
efficiently the number of lattice points on each $\gamma_n$, we
distinguish two cases based on the application of Lemma \ref{latticepts}.

\subsection*{ Case 1: $ {  2k \frac{N_2}{N_1} < \left(\frac{3}{4}\right)^{\frac{1}{3}}R_M^{-\frac{2}{3}} }  $ }

\

\noindent In this case \eqref{gamman}-\eqref{thetaM} guarantee
that the hypothesis of Lemma \ref{latticepts} is satisfied by each
arc of circle $\gamma_n$. Hence, on each $\gamma_n$ there are at
most two lattice points.

\subsection*{ Case 2: $ {  2k \frac{N_2}{N_1} \geq \left(\frac{3}{4}\right)^{\frac{1}{3}}R_M^{-\frac{2}{3}} }  $ }

\

\noindent In this case we approximate the number of lattice points
on $\gamma_n$  by the number of lattice points on
${\mathcal{C}}_{n}^{\lambda}$ (see for example \cite{bo1993} , \cite{bp}):
\begin{equation}\label{numblatt} 
\#{\mathcal{C}}_{n}^{\lambda} \lesssim R_M^{\epsilon} \sim (\lambda N_1)^{\epsilon} \lesssim (\lambda N_2)^{3 \epsilon}
\end{equation} 
for any $\epsilon > 0$.

%%%%%%%%%%%%%%%%%%%%%%%%%%%% COMMENT
{\comment{

\noindent \textbf{2a}: $N_1^2 \geq c.$ In this case $R_M \lesssim
\lambda N_1$. Therefore, from the assumption of Case 2, we
immediately get that $N_1 \lesssim \lambda^2 N_2^3$. Hence, $R_M
\lesssim (\lambda N_2)^3$ which combined with \eqref{numblatt}
gives:
\begin{equation}
\#{\gamma_n \lesssim (\lambda N_2)^{3 \epsilon}}\end{equation} for
any $\epsilon > 0$

\

\noindent \textbf{2b}: $N_1^2 < c.$ In this case we infer that
$R_M \lesssim \lambda$, which together with \eqref{numblatt}
implies:
\begin{equation}
\#{\gamma_n \lesssim \lambda^{\epsilon}}\end{equation} for any
$\epsilon > 0$

}}
%%%%%%%%%%%%%%%%%%%%%%%%%%%%%%%%%%%%%%% END OF COMMENT

\

\noindent Combining the estimate in \eqref{latpts}, {\bf{Case 1}}
and {\bf{Case 2}} we conclude that
$$\#S \lesssim \lambda^2 + \lambda^2 \; (\lambda N_2)^{\epsilon},$$
for any $\epsilon > 0$. Since $\lambda, N_2 \geq 1 $, together
with \eqref{2DM}, this implies that
$$M \lesssim (\lambda N_2)^{\epsilon},$$ for all
positive $\epsilon$'s. Hence \eqref{bila} follows.

\

(b) As in part (a) we will count the number of points in the set
$A^{\lambda}$ given by \eqref{Alam}. According to \eqref{B1} we
have that
\begin{equation} \label{scA}
A^{\lambda} = {\mathbb{Z}}^2 \cap \lambda B_{1},
\end{equation}
with
\begin{equation} \label{B1nol}
B_{1} = \{ (x,y) \in {\mathbb{R}}^2 : \; \left|x^2 + y^2 + 2(ax +
by)\right| \leq c,\; x^2 + y^2 \leq (k_2N_2)^{2},\; a^2 + b^2 \sim
N_1^2 \}.
\end{equation}

In the rest of the proof for arbitrary two points
$P_1, P_2$ on a circle with center $C$ we denote by 
$S(C,P_1,P_2)$ the solid sector contained between the lines
$CP_1$,$CP_2$ and the circle arc between $P_1$ and $ P_2$.

% Fig lam1
\begin{figure}[htbp]
\centering
\includegraphics*[scale=0.7,trim=0in 5.8in 3.5in 0in]{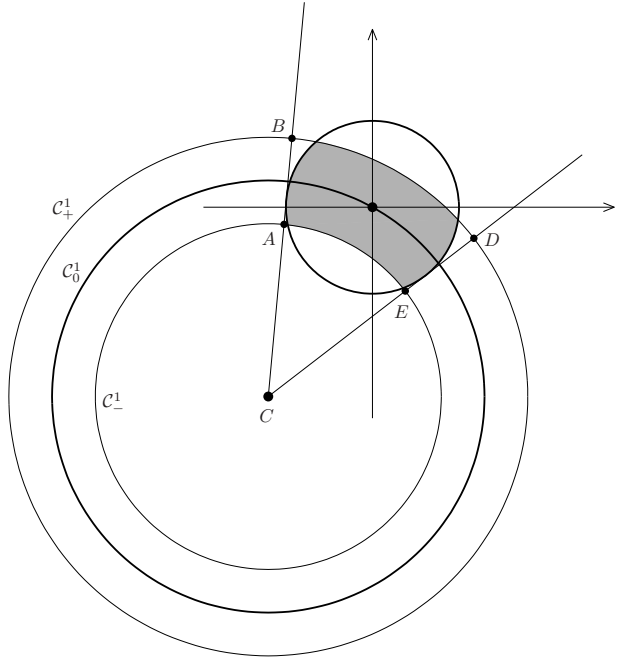}
\caption{Circular sectors when $\lambda = 1$}
\label{fig:5}
\end{figure}

We observe that (see Figure \ref{fig:5} )
$$B_1 \subset S(C,B,D) \setminus S(C,A,E),$$ where $CB$ and $CD$
lie along the tangent lines from $C=(-a,-b)$ to the circle $x^2 +
y^2 = (k_2N_2)^{2}$. Hence, setting $S_1=S(C,B,D), S_2=S(C,A,E)$,
and $H=S_1 \setminus S_2$ we have
\begin{align}
\# A^{\lambda} & = \# \{ {\mathbb{Z}}^{2} \cap \lambda B_1 \} \nonumber \\
& \leq \# \{ {\mathbb{Z}}^2 \cap \lambda S_1 \}
- \# \{ {\mathbb{Z}}^2 \cap \lambda S_2 \} \nonumber \\
& = \lambda^2 |S_1| - \lambda^2 |S_2| + O(\lambda) \label{useGauss} \\
& = \lambda^2 |H| + O(\lambda), \label{Alamcount}
\end{align}
where to obtain \eqref{useGauss} we used Lemma \ref{Gauss} twice.
In order to compute the area of $H$ notice that since $a^2 + b^2
\sim N_1^2$, we have (see Figure \ref{fig:5})
$$\sin\theta \sim \frac{N_2}{N_1}.$$ 
Hence, by the assumption $N_1 >> N_2$, we
get $\theta \sim \frac{N_2}{N_1}$. Therefore 
$$|H|= (r_1^2-r_2^2)\theta \sim \frac{N_2}{N_1},$$ 
with 
$$r_1^2 = c + (a^2 + b^2), \ \ r_2^2 = -c + (a^2+b^2).$$ 
This combined with
\eqref{Alamcount} implies that
$$\# A^{\lambda} \lesssim \lambda + \lambda^{2} \frac{N_2}{N_1},$$
that together with \eqref{2DM} gives
$$M \lesssim \left( \frac{1}{ {\lambda} } + \frac{N_2}{N_1} \right)^{1/2}.$$

\end{proof}

Now we are ready to prove the desired decay of the first modified
energy, which is the content of the following proposition:

\begin{prop}\label{energybound}Let $s>1/2, \; \lambda \lesssim N, \; t>0$, and $u_0 \in
H^s(\Bbb T_\lambda^2)$ be given. If  $u$ is a solution to
$(\ref{lambdaivp})-\eqref{lambdabc}$ then
$$|E^1(u)(t) - E^1(u)(0)| \lesssim \frac{1}{N^{1-}}\|Iu\|_{X^{1,1/2+}}.$$
\end{prop}

\begin{proof}The definition of $E$ together with equation \eqref{lambdaivp} give that:
\begin{align}\nonumber
\partial_t E(Iu)(t) &=
\textrm{Re}\int_{\Bbb T_\lambda^2}{\overline{I(u)_t}(|Iu|^2Iu -
\Delta Iu - iIu_t)}\\ \nonumber 
&= \textrm{Re}\int_{\Bbb T_\lambda ^2}{\overline{I(u)_t}(|Iu|^2 Iu - I(|u|^2 u))}
\end{align}
Therefore, integrating in time, and using the Parseval's formula
we get:
\begin{align}
&E^1(u)(t) - E^1(u)(0)= \\ \nonumber
&=\int_{0}^{t}\int_{\Gamma_4}\left(1-\frac{m(k_2+k_3+k_4)}{m(k_2)m(k_3)m(k_4)}\right)
\widehat{\overline{I\partial_tu}}(k_1)
\widehat{Iu}(k_2)\widehat{\overline{Iu}}(k_3)\widehat{Iu}(k_4
)
\end{align}
Using the equation \eqref{lambdaivp}, we then get:
$$E^1(u)(t) - E^1(u)(0)= \mbox{Tr}_1 + \mbox{Tr}_2,$$ where
\begin{align}
&\mbox{Tr}_1= \\ \nonumber
&=\int_{0}^{t}\int_{\Gamma_4}\left(1-\frac{m(k_2+k_3+k_4)}{m(k_2)m(k_3)m(k_4)}\right)
\widehat{\Delta \overline{Iu}}(k_1)\widehat{Iu}(k_2)\widehat{\overline{Iu}}(k_3)\widehat{Iu}(k_4
)
\end{align}
\begin{align}
&\mbox{Tr}_2= \\ \nonumber
&=\int_{0}^{t}\int_{\Gamma_4}\left(1-\frac{m(k_2+k_3+k_4)}{m(k_2)m(k_3)m(k_4)}\right)
\widehat{\overline{I
(|u|^2u)}}(k_1)\widehat{Iu}(k_2)\widehat{\overline{Iu}}(k_3)\widehat{Iu}(k_4
)
\end{align}
We wish to prove that:
\begin{equation}\label{Tr1Tr2}
|\mbox{Tr}_1| + |\mbox{Tr}_2| \lesssim \frac{1}{N^{1-}},
\end{equation}
for a constant $C=C(\|Iu\|_{X^{1,1/2+}}).$ According to the
\textbf{Decomposition Remark} in our estimates we may ignore the presence of the
complex conjugates.
\\
In order to obtain the estimate (\ref{Tr1Tr2}), we break $u$ into
a sum of dyadic pieces $u_j$, each with frequency support
$\langle k \rangle \sim 2^j$, $j=0,...$. We will then be able to
sum over all the frequency pieces $u_j$ of $u$, since our
estimates will decay geometrically in these frequencies.
\\
We start by analyzing $\mbox{Tr}_1$. First notice that,
$$\|\Delta (Iu)\|_{X^{-1,1/2+}} \leq \|Iu\|_{X^{1,1/2+}},$$
therefore, the estimate on $\mbox{Tr}_1$ will follow, if we prove that:
\begin{align}\label{term1bd}
&\left|\int_{0}^{t}
\int_{\Gamma_4}\left(1-\frac{m(k_2+k_3+k_4)}{m(k_2)m(k_3)m(k_4)}\right)
\widehat{\overline{\phi_1}}(k_1)\widehat{\phi_2}(k_2)\widehat{\overline{\phi_3}}(k_3)\widehat{\phi_4}(k_4)
\right| \\\nonumber
&\lesssim
\frac{1}{N^{1-}}
{(N_1N_2 N_3 N_4) }^{0-}
\|\phi_1\|_{X^{-1,1/2^+}}
\prod_{i=2}^{i=4}\|\phi_i\|_{X^{1,1/2+}}
\end{align}
for any $\lambda$-periodic function 
$\phi_i, i=1,...,4$ with
positive spatial Fourier transforms supported on
\begin{equation}
\langle k \rangle \sim 2^{l_i} \equiv N_i,
\end{equation}
for some $l_i \in \{0,1,...\}.$
\\
By the symmetry of the multiplier, we can assume that
\begin{equation}
N_2 \geq N_3 \geq N_4.
\end{equation}
Moreover, since we are integrating over $\Gamma_4$, we have that
$N_1 \lesssim N_2$. Therefore, we only need to get the factor
$N_2^{0-}$ in the estimate (\ref{term1bd}), in order to be able to
sum over all dyadic pieces. \\
Now we consider three different cases, by comparing $N$ to the size of the $N_i$'s.

\

\noindent\textbf{Case I.} $N >> N_2$. In this case, the symbol in
$\mbox{Tr}_1$ is identically zero, and the desired bound follows
trivially.

\

\noindent\textbf{Case II.} $N_2 \gtrsim N >> N_3 \geq N_4.$ Since
we are on $\Gamma_4,$ we also have that $N_1 \sim N_2.$ By the
mean value theorem, we have
$$\left|1- \frac{m(k_2+k_3+k_4)}{m(k_2)}\right|
= \left| \frac{m(k_2) - m(k_2 + k_3 + k_4)}{m(k_2)} \right|
\lesssim \left| \frac{m^{\prime}(k_2)}{m(k_2)} \right| N_3
\sim \frac{N_3}{N_2}.$$
Using the pointwise bound above, the Cauchy-Schwartz inequality,
and Plancharel's theorem we get:
\begin{equation}
|\mbox{Tr}_1| \lesssim \frac{N_3}{N_2}\|\phi_1
\phi_3\|_{L^2_tL^2_x}\|\phi_2 \phi_4\|_{L^2_tL^2_x}
\end{equation}
By two applications of H\"{o}lder's inequality, followed by four
applications of the Strichartz estimate (\ref{L44}) we obtain:
$$|\mbox{Tr}_1| \lesssim
\frac{N_3}{N_2} {\lambda}^{0+} \frac{(N_1)^{1+}}{(N_2N_3N_4)^{1-}}
\|\phi_1\|_{X^{-1,1/2+}}
\|\phi_3\|_{X^{1,1/2+}}\|\phi_2\|_{X^{1,1/2+}}
\|\phi_4\|_{X^{1,1/2+}}.$$ Since $N_1 \sim N_2 \gtrsim N$, and
$\lambda \lesssim N$, it follows immediately that in this case
$$|\mbox{Tr}_1| \lesssim \frac{1}{N^{1-}}N_2^{0-}\|\phi_1\|_{X^{-1,1/2+}}
\|\phi_2\|_{X^{1,1/2+}}\|\phi_3\|_{X^{1,1/2+}}
\|\phi_4\|_{X^{1,1/2+}}.$$
\\
\textbf{Case III.} $N_2 \geq N_3 \gtrsim N.$ We will use the
following crude bound on the multiplier:
\begin{equation}\label{mult}
\left|1-\frac{m(k_2+k_3+k_4)}{m(k_2)m(k_3)m(k_4)}\right|
\lesssim \frac{m(k_1)}{m(k_2)m(k_3)m(k_4)},\end{equation}
which follows from the fact that we are integrating on $\Gamma_4$.
We also need the following estimate:
\begin{equation}\label{monot} \frac{1}{m(k)|k|^{\alpha}}
\lesssim N^{-\alpha},
\end{equation}
for $\alpha \geq 1/2-$, $|k| \gtrsim N$.
\\
We distinguish two subcases:

\

\noindent$\textbf{IIIa)} N_2 \sim N_3 \gtrsim N.$ Using
Plancharel's, followed by Holder's inequality, and then by four
applications of the Strichartz estimate (\ref{L44}), we get:
\begin{align}
|\mbox{Tr}_1| & \lesssim \frac{m(k_1)}{m(k_2)m(k_3)m(k_4)}
{\lambda}^{0+}
\frac{(N_1)^{1+}}{(N_2N_3N_4)^{1-}} \times \nonumber \\
& \times \|\phi_1\|_{X^{-1,1/2+}}\|\phi_2\|_{X^{1,1/2+}}
\|\phi_3\|_{X^{1,1/2+}} \|\phi_4\|_{X^{1,1/2+}}\nonumber \\
& \lesssim \lambda^{0+} \frac{1}{N^{1-} m(k_4) N_4^{1-}}
\|\phi_1\|_{X^{-1,1/2+}} \|\phi_2\|_{X^{1,1/2+}}
\|\phi_3\|_{X^{1,1/2+}} \|\phi_4\|_{X^{1,1/2+}},\label{usemon}
\end{align}
where to obtain \eqref{usemon} we use that $N_1 \leq N_2$ together
with \eqref{monot} with $\alpha = 1-$ ($N_3 \gtrsim N$). Now if
$N_4 \gtrsim N$, we conclude with an application of the
monotonicity property (\ref{monot}). If $N_4 << N$, then we use
that $m(k_4)=1$ and $N_4 \geq 1.$ Hence, since $\lambda \lesssim
N$, we get
$$|\mbox{Tr}_1| \lesssim
\frac{1}{N^{1-}}N_2^{0-} \|\phi_1\|_{X^{-1,1/2+}}
\|\phi_2\|_{X^{1,1/2+}} \|\phi_3\|_{X^{1,1/2+}}
\|\phi_4\|_{X^{1,1/2+}}$$

\noindent$\textbf{IIIb)} N_2>>N_3 \gtrsim N$. Since we are
integrating on $\Gamma_4$, we get $N_1 \sim N_2$. 
Let us emphasize that in the case \textbf{IIIb)} we could not use
the same strategy as in the case \textbf{IIIa)}. More precisely,
the approach of \textbf{IIIa)} consisting of
Plancharel's, followed by Holder's inequality, and then by four
applications of the Strichartz estimate (\ref{L44})
would produce a bound of the type
\begin{align*}
& \frac{m(k_1)}{m(k_2)m(k_3)m(k_4)} {\lambda}^{0+}
\frac{(N_1)^{1+}}{(N_2N_3N_4)^{1-}}
\times \|\phi_1\|_{X^{-1,1/2+}}\|\phi_2\|_{X^{1,1/2+}}
\|\phi_3\|_{X^{1,1/2+}} \|\phi_4\|_{X^{1,1/2+}} \\
& \lesssim \frac{1}{N^{1-}}\lambda^{0+} N_1^{0+}
\|\phi_1\|_{X^{-1,1/2+}} \|\phi_2\|_{X^{1,1/2+}}
\|\phi_3\|_{X^{1,1/2+}} \|\phi_4\|_{X^{1,1/2+}}
\end{align*}
which we cannot sum in $N_1$ (here we cannot cancel $N_1^{0+}$ by $N_3^{0 +}$ since
$N_3 << N_1$). Thus to handle the
case \textbf{IIIb)} we do need an improved Strichartz estimate
given in Proposition \ref{2Dcounting}.

In order to obtain the desired estimate, we apply the
Cauchy-Schwartz inequality, together with Proposition
\ref{2Dcounting} for the pair $\phi_1,\phi_3$ and for the pair
$\phi_2,\phi_4$. Thus
\begin{align} |\mbox{Tr}_1| &\lesssim
\frac{m(k_1)}{m(k_2)m(k_3)m(k_4)}(\lambda N_3)^{\epsilon} (\lambda N_4)^{\epsilon}
\frac{N_1}{N_2N_3N_4} \times \nonumber\\
&\times \|\phi_1\|_{X^{-1,1/2+}}\|\phi_2\|_{X^{1,1/2+}}
\|\phi_3\|_{X^{1,1/2+}} \|\phi_4\|_{X^{1,1/2+}} \lesssim \nonumber\\
&\lesssim
\frac{1}{m(k_3)m(k_4)} \; \lambda^{0+} \;
\frac{1}{N_3^{1-} N_4^{1-}} \times \label{IIIbusesmall2}\\
&\times \|\phi_1\|_{X^{-1,1/2+}}\|\phi_2\|_{X^{1,1/2+}}\|\phi_3\|_{X^{1,1/2+}} \|\phi_4\|_{X^{1,1/2+}} \lesssim \nonumber\\
&\lesssim
\frac{1}{N^{1-}} \lambda^{0+} \frac{1}{m(k_4)N_4^{1-}}
\|\phi_1\|_{X^{-1,1/2+}}\|\phi_2\|_{X^{1,1/2+}}\|\phi_3\|_{X^{1,1/2+}}
\|\phi_4\|_{X^{1,1/2+}}\label{IIIbusemon3}
\end{align}
where in order to obtain \eqref{IIIbusesmall2} we have used that
$N_1 \sim N_2$ and in order to obtain
\eqref{IIIbusemon3} we have used \eqref{monot} once with $\alpha = 1-$ ($N_3 \gtrsim N$).
Finally  as observed in the case $\textbf{IIIa)}$, $m(k_4) N_4^{1-}\geq 1$. Therefore
$$ |\mbox{Tr}_1| \lesssim
\frac{1}{ N^{1-} } N_2^{0-}
\|\phi_1\|_{X^{-1,1/2+}}\|\phi_2\|_{X^{1,1/2+}}
\|\phi_3\|_{X^{1,1/2+}} \|\phi_4\|_{X^{1,1/2+}}.$$

%%%%%%%%%%%%%%%%%%%%%%%%%%%%%%%%%%%%%%%%%%% COMMENT

\comment{In order to obtain the desired estimate, we apply the
Cauchy-Schwartz inequality, together with Proposition
\ref{2Dcounting} for the pair $\phi_1,\phi_3$ and for the pair
$\phi_2,\phi_4$. Thus
\begin{align} |\mbox{Tr}_1| &\lesssim
\frac{m(k_1)}{m(k_2)m(k_3)m(k_4)}(\frac{1}{\lambda^2}+\frac{1}{N_2})
\frac{N_1}{N_2N_3N_4} \times \nonumber\\
&\times \|\phi_1\|_{X^{-1,1/2+}}\|\phi_2\|_{X^{1,1/2+}}
\|\phi_3\|_{X^{1,1/2+}} \|\phi_4\|_{X^{1,1/2+}} \lesssim \nonumber\\
&\lesssim
\frac{1}{m(k_3)m(k_4)}(\frac{1}{\lambda^2}+\frac{1}{N})
\frac{1}{N_3 N_4} \times \label{IIIbusesmall2}\\
&\times \|\phi_1\|_{X^{-1,1/2+}}\|\phi_2\|_{X^{1,1/2+}}\|\phi_3\|_{X^{1,1/2+}} \|\phi_4\|_{X^{1,1/2+}} \lesssim \nonumber\\
&\lesssim
\frac{1}{N}(\frac{1}{\lambda^2}+\frac{1}{N})\frac{1}{m(k_4)N_4}
\|\phi_1\|_{X^{-1,1/2+}}\|\phi_2\|_{X^{1,1/2+}}\|\phi_3\|_{X^{1,1/2+}}
\|\phi_4\|_{X^{1,1/2+}}\label{IIIbusemon3}
\end{align}
where in order to obtain \eqref{IIIbusesmall2} we have used that
$N_1 \sim N_2 \gtrsim N$ and in order to obtain
\eqref{IIIbusemon3} we have used \eqref{monot} once.
Finally  as observed in the case $\textbf{a1.2)}$, $m(k_4) N_4
\geq 1$. Therefore
$$ |\mbox{Tr}_1| \lesssim
\frac{1}{N}(\frac{1}{\lambda^2}+\frac{1}{N})
\|\phi_1\|_{X^{-1,1/2+}}\|\phi_2\|_{X^{1,1/2+}}
\|\phi_3\|_{X^{1,1/2+}} \|\phi_4\|_{X^{1,1/2+}}.$$}

%%%%%%%%%%%%%%%%%%%%%%%%%%%%%%%%%%%%%%%%%%%%%%%%%%% END OF COMMENT

Now, we
proceed to analyze $\mbox{Tr}_2$. We claim that the following
stronger estimate holds:
$$|\mbox{Tr}_2| \lesssim \frac{1}{N^{2-}}\|Iu\|^6_{X^{1,1/2+}}.$$
As for $\mbox{Tr}_1$ it suffices to prove the following bound:
$$\mbox{LHS}=\int_0^{t}\int_{\Gamma_6}m_{123}(m_{456}- m_4m_5m_6)\prod_{j=1}^{6} \widehat{\phi_j}
\lesssim
\frac{N_{max}^{0-}}{N^{2-}}\prod_{i=1}^{6}\|I\phi_i\|_{X^{1,1/2+}},$$
with the $\phi_i$'s having positive spatial Fourier transform,
supported on $$\langle k \rangle \sim 2^{l_i}\equiv N_i,$$ for
some $l_i \in \{0,1,...\}$, and with $N_{max}, N_{med}$ denoting
respectively the biggest and the second biggest frequency among
the $N_i$'s.\\
Since we are integrating over $\Gamma_6$, we can assume that
$N_{max} \sim N_{med}$. Also, we only need to consider the case
when $N_{max} \gtrsim N$, otherwise the symbol in (LHS)
is zero, and the estimate above holds trivially. \\
Now, using Holder's inequality, together with the estimate
(\ref{monot}) twice, and the fact that the multiplier $m$ is bounded, we
then obtain:
\begin{align*}
\mbox{LHS} & \leq \frac{N_{max}^{1-}N_{med}^{1-}m_{max}m_{med}}
{N_{max}^{1-}N_{med}^{1-}m_{max}m_{med}}\int_0^{t}\int_{\Gamma_6}\prod_{j=1}^{6}
\widehat{\phi_j}(k,t) \\
& \lesssim
\frac{N_{max}^{0-}}{N^{2-}}\|J^{1-}I\phi_{max}\|_{L^4_tL^4_x}
\|J^{1-}I\phi_{med}\|_{L^4_tL^4_x}
\prod_{k=1}^{4}\|\phi_k\|_{L^8_tL^8_x}.
\end{align*}
The Strichartz estimate
(\ref{L44}), and the interpolation estimate (\ref{interpol}) then imply:
$$\mbox{LHS} \lesssim  \frac{N_{max}^{0-}}{N^{2-}}\|I\phi_{max}\|_{X^{1,1/2+}}\|I\phi_{med}\|_{X^{1,1/2+}}
\prod_{k=1}^{4}\|\phi_k\|_{X^{1/2+,1/2+}} \lesssim
 \frac{N_{max}^{0-}}{N^{2-}} \prod_{j=1}^{6}\|I\phi_j\|_{X^{1,1/2+}},$$ where in the last
inequality we have used (\ref{Iprop}).
\end{proof}

\subsection*{Step 3: Proof of Theorem \ref{main2D}} We are now ready to 
exhibit the proof of Theorem \ref{main2D}. As we explained in 
Theorem \ref{main1D} we do not keep track of the $\lambda^{0+}$ factors.

\begin{proof} Let $u_0$ be an initial data, and let $N>>1$ be
given. Recall that, for any $\lambda \in \Bbb R$, we set
$u_0^\lambda(x,t)=
\frac{1}{\lambda}u_0(\frac{x}{\lambda},\frac{t}{\lambda^2}).$ Now,
let us choose a  rescaling parameter $\lambda$ so that
$E^1(u^\lambda_0) \lesssim 1$, that is $$\lambda \sim
N^{\frac{1-s}{s}}.$$ Then, by Corollary \ref{n2Dlwp} and Proposition \ref{energybound},
there exists a $\delta$, and a solution $u_\lambda$ to
(\ref{lambdaivp}), with initial condition $u_0^{\lambda}$, such
that
\begin{align*}
E^1(u_\lambda)(\delta) & \lesssim E^1(u_\lambda)(0) +
O\left(\frac{1}{N^{1-}}\right).
\end{align*}
Therefore, we can continue our solution $u_\lambda$ until the size
of $E^1(u_\lambda)(t)$ reaches 1, that is at least $C \cdot
N^{1-}$ times. Hence
$$E^1(u_\lambda)(C \cdot N^{1-} \delta) \sim 1.$$ Given
$T>>1$, we choose our parameter $N>>1$, so that $$T \sim \frac{C
\cdot N^{1-} \delta}{\lambda^2} \sim N^{\frac{3s-2}{s}-}.$$ We
observe that the exponent of $N$ is positive as long as $s>2/3$,
hence $N$ is well-defined for all times $T.$ Now, undoing the
scaling, we get that:
$$E^1(u)(T) \lesssim T^{\frac{4(1-s)}{3s-2}+},$$ with $u$ solution to (\ref{ivp1})-(\ref{bc1}), $d=2$.
This bound implies the desired conclusion.
\end{proof}

\end{document}